\documentclass[a4paper,10pt]{amsart}
\usepackage[german,english]{babel}
\usepackage[latin1]{inputenc}
\usepackage{ucs}
\usepackage{amsmath}
\usepackage{amsfonts}
\usepackage{amssymb}
\usepackage{amsthm}
\usepackage{fontenc}
\usepackage{graphicx}
\usepackage{url}
\usepackage{color}
\usepackage[T1]{fontenc}
\usepackage{a4wide}

\newtheorem{theorem}{Theorem}

\newtheorem{proposition}{Proposition}

\newtheorem{remark}{Remark}

\newcommand{\bsg}{\boldsymbol{g}}
\newcommand {\intQT} {\int_{Q}} 
\newcommand {\intOm}  {\int_\Omega}

\def \dx     {d\boldsymbol{x}}
\def \dxdt    {d\boldsymbol{x}\,dt}

\title[Guaranteed lower bounds for cost functionals]
{Guaranteed lower bounds for cost functionals of time-periodic parabolic 
 optimization problems}

\author{Monika Wolfmayr}

\thanks{This research was supported by the Academy of Finland, grant 295897.}

\address[M. Wolfmayr]{Faculty of Information Technology; University of Jyv\"askyl\"a; 
P.O. Box 35 (Agora); FIN-40014 Jyv\"askyl\"a; Finland}
\email{monika.k.wolfmayr@jyu.fi}

\begin{document}

\begin{abstract}
In this paper, a new technique is
shown for 
deriving computable, guaranteed lower bounds of functional type (minorants)
for two different cost functionals subject to a parabolic time-periodic boundary value problem.
Together with previous results on upper bounds (majorants) for one of the cost functionals,
both minorants and majorants lead to
two-sided estimates of functional type 
for the optimal control problem. 
Both upper and lower bounds are
derived for the second new cost functional
subject to 
the same parabolic 
PDE-constraints, but where the target is a desired gradient.
The time-periodic optimal control problems are discretized by the multiharmonic finite element method
leading to large systems of linear equations  having a saddle point structure.
The derivation of preconditioners for the minimal residual method for the new optimization
problem is discussed in more detail.
Finally, several numerical experiments for both optimal control problems
are presented confirming the theoretical results obtained.
This work provides the basis for an adaptive scheme for time-periodic optimization problems.
\end{abstract}

\maketitle

\section{Introduction}
\label{Sec1:Intro}

In this work, we derive fully computable, guaranteed lower bounds (minorants)
for cost functionals of parabolic optimization 
problems with given time-periodic conditions. 
The fully computable upper bounds (majorants) for one of the cost functionals were derived in 
\cite{LangerRepinWolfmayr2016}. The second cost functional
is new and for this one both upper and lower bounds are presented.
The motivation for the second problem lies in applications, where the target is the gradient or 
flux instead of the state function.
Lower bounds for cost functionals
of time-periodic parabolic optimal control problems have not been discussed yet.
However, optimal control problems are highly important for different applications
(see e.g. \cite{NeittaanmaekiSprekelsTiba:2006} as well as the original work
\cite{Lions:1971}). 
These applications also include time-periodic problems, see for instance
\cite{AltmannStingelinTroeltzsch:2014} and \cite{HouskaLogistImpeDiehl:2009}
considering problems in electromagnetics and biochemistry, respectively.
For these types of problems the multiharmonic (or harmonic-balanced) finite element method
(short MhFEM) is a natural choice. The functions are approximated by truncated Fourier series
in time and by the finite element method (FEM) in space -- more precisely, the Fourier coefficients.
We refer to the application of this discretization technique  
already in \cite{YamadaBessho:1988} as well as later in 
\cite{BachingerKaltenbacherReitzinger:2002a,
BachingerLangerSchoeberl:2005,
BachingerLangerSchoeberl:2006,
CopelandLanger:2010} 
for non-linear time-harmonic eddy current problems.
Moreover, time-periodic optimal control problems and the MhFEM were discussed in e.g.
\cite{KollmannKolmbauer:2011,
KollmannKolmbauerLangerWolfmayrZulehner:2013, 
LangerWolfmayr:2013} 
and 
\cite{KolmbauerLanger:2012, KolmbauerLanger:2013}.
Recent works on robust preconditioners for time-periodic parabolic and eddy-current optimal
control problems are discussed in
\cite{Liang2018} and \cite{AxelssonLukas2018}, respectively.
In this work, a standard finite element discretization is used with continuous, piecewise linear
finite elements and a regular grid as discussed, e.g., in
\cite{Braess:2005, Ciarlet:1978, Steinbach:2008}. 
However, the method is wider applicable using also other finite elements or also finite differences
(for instances, if the given domain is geometrically rather non-complex).

Functional a posteriori estimation 
provides a useful machinery 
to derive computable and guaranteed
quantities for the desired unknown solution, see, e.g.,
\cite{Repin2002,
GaevskayaRepin:2005} on parabolic problems.
Recent works on new estimates for parabolic problems and parabolic optimal control problems
can be found in \cite{Matculevich2016} and \cite{Shakya2017}, respectively.
A posteriori estimates of functional type 
for elliptic optimal control problems can be found in
\cite{GaevskayaHoppeRepin:2006, GaevskayaHoppeRepin:2007, Repin2008,
MaliNeittaanmaekiRepin:2014}.
First functional type estimates for inverse problems, 
which are related to
optimal control problems, can be found 
in \cite{RepinRossi2013, ClasonKaltenbacherWachsmuth2016}.
Moreover, recent results on guaranteed computable estimates for convection-dominated diffusion problems
are presented in \cite{MatWolf2018}.
In \cite{LangerRepinWolfmayr2016}, majorants 
for one cost functional of a time-periodic parabolic optimal control and for the 
corresponding optimality system were presented.
This work presents
the corresponding minorants
for this cost functional
using the new technique presented in  
\cite{Wolfmayr2016}, which makes use of ideas  
derived by Mikhlin \cite{Mikhlin1964}
but generalized for the class of optimal control problems. 
We mention here that
\cite{Mali2017} presents a different approach for the derivation
of a lower bound for a class of elliptic optimal control problems.

We extend the analysis in this paper and consider a second cost functional with respect
to the same parabolic time-periodic boundary value problem.
In the second optimal control problem, the target is a given desired gradient. Problems of that type
have been earlier discussed in \cite{GaevskayaHoppeRepin:2007}.
The results on computable lower bounds together with the upper bounds
lead to two-sided estimates which can be used to derive majorants for the discretization
error in state and control. These majorants and minorants provide a new formulation
of the optimization problems since they can, in principle, be used as objects
of direct minimization on their difference.
The majorants and minorants can be used in order to derive an adaptive scheme in time and in space.
The linearity of the optimal control problems leads to decoupled problems in the Fourier modes
including decoupling in the majorants and minorants, which is shown in this work.
The overall estimators provide the modes/mode numbers
which are necessary for computations. The problems for the different modes can then be
computed on different grids, for which estimators in space can be used with respect to
the finite element discretized Fourier coefficients. Altogether we derive a space-time adaptive method.
Its idea has been for the first time introduced in \cite{LangerRepinWolfmayr2016} and has been called
the adaptive multiharmonic finite element method (AMhFEM).

In this work, robust preconditioners for the preconditioned minimal residual (MinRes) method (see
\cite{PaigeSaunders:1975}) are discussed 
for the second optimization problem, which 
are new for this case. 
Also the practical performance of the AMLI preconditioner MinRes solver is
presented in various numerical experiments for both optimization problems.
For additional numerical tests regarding the AMLI precondititioned MinRes solver used in this work
and its performance for different cases of given data in time-periodic parabolic optimal control problems,
we refer the reader to \cite{LangerWolfmayr:2013}.

The work is arranged in the following sections:
In Section~\ref{Sec:ModelProblems}, the two types of cost functionals are presented 
including some preliminary results.
We denote the problems by optimization problem I and II.
Moreover, the former result on the majorant for problem I is summarized there.
The new minorant for optimization problem I is derived in
Section \ref{Sec:TwoSidedBoundsI} followed by the discussion of
the majorant and minorant for optimization problem II
in Section \ref{Sec:TwoSidedBoundsII}.
In Section \ref{Sec:MhFEM}, the MhFEM for both 
optimization problems is presented.
Robust preconditioners
for applying the preconditioned MinRes method on the problems discretized by the MhFEM
are presented in Section \ref{Sec:Preconditioning}.
Section~\ref{Sec:NumericalResults} discusses detailed
a set of various numerical experiments for both optimization problems
I and II. A few final remarks are drawn in 
Section~\ref{Sec7:Conclusions}.

\section{Time-periodic parabolic PDE, the two cost functionals and preliminary results}
\label{Sec:ModelProblems}

We denote by $\Omega \subset \mathbb{R}^d$ with possible dimensions 
$d=\{1,2,3\}$ the spatial 
bounded Lipschitz domain with boundary $\Gamma := \partial \Omega$.
Also we denote by $Q := \Omega \times (0,T)$ the space-time domain
and $\Sigma := \Gamma \times (0,T)$ its lateral surface, 
where $(0,T)$ is the given time interval.
The optimization problems are both subject to the following 
parabolic PDE with given homogeneous Dirichlet boundary conditions and
time-periodical condition: 
\begin{align}
\label{equation:forwardpde:OCP1}
 \sigma \, \partial_t y
 - \nabla \cdot ( \nu \nabla y)
 &= u 
 \hspace{0.5cm} &\text{in }  
 Q, \\
\label{equation:forwardpde:OCP2}
 y
 &= 0 \hspace{0.5cm} &\text{on } 
 \Sigma, \\
\label{equation:forwardpde:OCP3}
 y(0) &= y(T) \hspace{0.5cm} &\text{in } 
 \overline{\Omega}.
\end{align}
The function $y$ is the state and $u$ will be the control function.
We assume that the coefficient functions 
$\sigma$ and $\nu$
are positive and  
bounded satisfying 
 $0< \underline{\sigma} \leq \sigma(\boldsymbol{x}) \leq \overline{\sigma}$
and
 $0< \underline{\nu} \leq \nu(\boldsymbol{x}) \leq \overline{\nu}$
for $\boldsymbol{x} \in \Omega$
with constants
$\underline{\sigma}$, $\overline{\sigma}$, $\underline{\nu}$ and $\overline{\nu}$.
The time-periodic problems in this paper are motived by real-life applications such as
computational electromagnetics, where these parameters correspond to the reluctivity and conductivity
being usually piecewise constant because of the various materials of the electrical devices.

\subsection{Preliminaries}
\label{SSec:Preliminaries}

In the following, we present a
proper functional space setting for time-periodic problems which starts by defining
the Hilbert spaces
\begin{align*}
 H^{1,0}(Q) &= \{u \in L^2(Q) : \nabla u \in [L^2(Q)]^d \}, \,
 H^{1,0}_0(Q) = \{u \in H^{1,0}(Q) : u = 0 \text{ on } \Sigma\}, \\
 H^{0,1}(Q) &= \{u \in L^2(Q) : \partial_t u \in L^2(Q) \}, \,
 H^{0,1}_{per}(Q) = \{u \in H^{0,1}(Q) :   u(0) = u(T) \text{ in } \overline{\Omega} \}, \\
 H^{1,1}(Q) &= \{u \in L^2(Q) : \nabla u \in [L^2(Q)]^d, \partial_t u \in L^2(Q) \},
\end{align*}
(see also \cite{Ladyzhenskaya:1973, LadyzhenskayaSolonnikovUralceva:1968}).
For instance, the norm in $H^{1,1}$ is given by
\begin{align*}
 \|u\|_{1,1} &:= \left(\int_{Q} \left(u(\boldsymbol{x},t)^2 + |\nabla u(\boldsymbol{x},t)|^2
			+ |\partial_t u(\boldsymbol{x},t)|^2 \right) \, d\boldsymbol{x} \, dt \right)^{1/2}.
\end{align*}
In the following, we skip the subindex for the norms and
inner products 
in $L^2(Q)$ 
and denote them by $\| \cdot \|$ and
$\langle \cdot, \cdot \rangle$.
For $L^2(\Omega)$ and $H^1(\Omega)$, we denote them by
$\| \cdot \|_{\Omega}$ and $\langle \cdot, \cdot \rangle_{\Omega}$
and $\| \cdot \|_{1,\Omega}$ and
$\langle \cdot, \cdot \rangle_{1,\Omega}$,
respectively.

Time-periodic functions which are at least from the space $L^2$ can be naturally represented
by Fourier series as
\begin{align*}
 u(\boldsymbol{x},t) = u_0^c(\boldsymbol{x}) +
      \sum_{k=1}^{\infty} \left(u_k^c(\boldsymbol{x}) \cos(k \omega t)
                          + u_k^s(\boldsymbol{x}) \sin(k \omega t)\right)
\end{align*}
for $\omega = 2 \pi /T$ being the frequency, $T$ the period and with the Fourier
cofficients
\begin{align*}
 \begin{aligned}
 u_0^c(\boldsymbol{x}) = \frac{1}{T} \int_0^T u(\boldsymbol{x},t) \,dt, \,
 u_k^c(\boldsymbol{x}) = \frac{2}{T} \int_0^T u(\boldsymbol{x},t) \cos(k \omega t)\,dt, \,
 u_k^s(\boldsymbol{x}) = \frac{2}{T} \int_0^T u(\boldsymbol{x},t) \sin(k \omega t)\,dt.
 \end{aligned}
\end{align*}
We define the norm in Fourier space as follows
\begin{align}
\label{definition:H01/2seminorm}
 \big\| \partial^{1/2}_t u \big\| ^2 :=
 |u|_{0,\frac{1}{2}}^2 := 
 \frac{T}{2} \sum_{k=1}^{\infty} k \omega \|\boldsymbol{u}_k\|_{\Omega}^2
\end{align}
as well as the spaces
$ H^{0,\frac{1}{2}}_{per}(Q) = \{ u \in L^2(Q) : \big\| \partial^{1/2}_t u \big\|  < \infty \}$, 
$ H^{1,\frac{1}{2}}_{per}(Q) = \{ u \in H^{1,0}(Q) : \big\| \partial^{1/2}_t u \big\|  < \infty \}$, 
$ H^{1,\frac{1}{2}}_{0,per}(Q) = \{ u \in H^{1,\frac{1}{2}}_{per}(Q): u = 0 \mbox{ on } \Sigma \}$,
where $\boldsymbol{u}_k = (u_k^c,u_k^s)^T$, $k \in \mathbb{N}$.
We also introduce the orthogonal vector
$\boldsymbol{u}_k^\perp = (-u_k^s,u_k^c)^T$.
The inner products (including also a $\sigma$-weighted version) in these spaces are defined by
 $\langle \partial^{1/2}_t u, \partial^{1/2}_t v \rangle  :=
 \frac{T}{2} \sum_{k=1}^{\infty} k \omega \langle\boldsymbol{u}_k,\boldsymbol{v}_k\rangle_{\Omega}$,
and
$ \langle\sigma
 \partial^{1/2}_t u, \partial^{1/2}_t v\rangle  :=
 \frac{T}{2} \sum_{k=1}^{\infty} k \omega \langle\sigma
 \boldsymbol{u}_k,\boldsymbol{v}_k\rangle_{\Omega}$.
The
$H^{1,\frac{1}{2}}_{per}(Q)$-seminorm is defined by
\begin{align*}
 |u|_{1,\frac{1}{2}}^2 
 = T \, \|\nabla u_0^c\|_{\Omega}^2
 + \frac{T}{2} \sum_{k=1}^{\infty} \left(k\omega \|\boldsymbol{u}_k\|_{\Omega}^2
 + \|\nabla \boldsymbol{u}_k\|_{\Omega}^2\right)
 = \|\nabla u\| ^2 + \|\partial_t^{1/2} u\| ^2 
\end{align*}
and the corresponding norm is
$\|u\|_{1,\frac{1}{2}}^2 := \|u\| ^2 + |u|_{1,\frac{1}{2}}^2$.
Finally, we also define the product
$ \langle \kappa,\partial_t^{1/2} u \rangle
 := \frac{T}{2} \sum_{k=1}^{\infty} (k \omega)^{1/2} \langle \boldsymbol{\kappa}_k,\boldsymbol{u}_k \rangle_{\Omega}$
as well as the orthogonal Fourier series
\begin{align*}
 \begin{aligned}
 u^{\perp}(\boldsymbol{x},t) &:= 
 \sum_{k=1}^{\infty} (-\boldsymbol{u}_k^{\perp})^T
 \cdot (\cos(k \omega t),\sin(k \omega t))^T.
 \end{aligned}
\end{align*}
Using this notation we can prove that
 $(u^\perp)^\perp$ = -u,
 $
 \big\|u^\perp\big\|  = \big\|u\big\|$
 and
$
 \big\|\partial^{1/2}_t u^\perp\big\| 
 = \big\|\partial^{1/2}_t u\big\|$
 for all  
 $u \in H^{0,\frac{1}{2}}_{per}(Q)$
and also that $\|\boldsymbol{u}_k^\perp\|_{\Omega} = \|\boldsymbol{u}_k\|_{\Omega}$.
Briefly, we recall from \cite{LangerRepinWolfmayr2015}, the following
identities
 \begin{align}
 \label{equation:H11/2identities}
 \begin{aligned}
  \langle \sigma \partial_t^{1/2} u,\partial_t^{1/2} v  \rangle  =
  \langle \sigma \partial_t u,v^{\perp} \rangle, \,
  \langle \sigma \partial_t^{1/2} u,\partial_t^{1/2} v^{\perp}  \rangle =
  \langle \sigma \partial_t u,v  \rangle, \quad \forall u \in H^{0,1}_{per}(Q) \, \forall
  v \in H^{0,\frac{1}{2}}_{per}(Q),
 \end{aligned}
 \end{align}
as well as 
orthogonality relations 
\begin{align}
\label{equation:orthorelation1}
\langle \sigma \partial_t u,u \rangle  = 0 \, &\mbox{ and } \,
 \langle \sigma u^{\perp},u \rangle  = 0 \qquad &&\forall \, u \in H^{0,1}_{per}(Q), \\
\label{equation:orthorelation2}
 \langle \sigma \partial^{1/2}_t u,\partial^{1/2}_t u^{\perp} \rangle  = 0 \, &\mbox{ and } \,
\langle \nu \nabla u, \nabla u^{\perp} \rangle  = 0 \qquad &&\forall \, u \in H^{1,\frac{1}{2}}_{per}(Q).
\end{align}
Note that the following identity is valid
(in Fourier series sense)
\begin{align}
\label{def:identityH01/2per}
\int_{Q}
 \kappa \, \partial_t^{1/2} u^\perp \, d\boldsymbol{x}\,dt = - \int_{Q}
 \partial_t^{1/2} \kappa^\perp \, u \, d\boldsymbol{x}\,dt
 \qquad \forall \, \kappa, u \in H^{0,\frac{1}{2}}_{per}(Q).
\end{align}
Friedrichs' inequality 
in $Q$ can be proved by using standard Friedrichs' inequality
on the Fourier coefficients
with respect to the spatial domain $\Omega$. 
We have that
$ \|\nabla u\| ^2 
 \geq 
  \frac{1}{C_F^2} \|u\| ^2$,
where $C_F > 0$ is Friedrichs' constant.

In the following, the parameter $\lambda>0$ denotes the regularization or cost parameter.

\subsection{Optimization problem I}
\label{SSec:ModelProblemI}

In the first case, we want to minimize the following cost functional with respect to the unknown
state $y$ and control $u$:
\begin{align}
\label{equation:minfunc1}
 \mathcal{J}(y,u) :=  \frac{1}{2} \| y - y_d \|^2 + \frac{\lambda}{2}  \| u \|^2
\end{align}
subject to the time-periodic boundary value problem \eqref{equation:forwardpde:OCP1}--\eqref{equation:forwardpde:OCP3}.
The given desired state $y_d \in L^2(Q)$
does not have to be time-periodic. It only has to be from the space
$L^2(Q)$.
The cost functional $\mathcal{J}$ defined in \eqref{equation:minfunc1} can be written as
 \begin{align*}
  \mathcal{J}(y,u) 
	=  T \mathcal{J}_0(y_0^c,u_0^c)
	+ \frac{T}{2} \sum_{k=1}^\infty \mathcal{J}_k(\boldsymbol{y}_k,\boldsymbol{u}_k),
 \end{align*}
 where $\mathcal{J}_0(y_0^c,u_0^c)
  = \frac{1}{2} \|y_0^c - {y_d^c}_0\|_{\Omega}^2
     + \frac{\lambda}{2} \|u_0^c\|_{\Omega}^2$
 and
 $\mathcal{J}_k(\boldsymbol{y}_k,\boldsymbol{u}_k)
 = \frac{1}{2} \|\boldsymbol{y}_k-{\boldsymbol{y}_d}_k\|_{\Omega}^2
     + \frac{\lambda}{2} \|\boldsymbol{u}_k\|_{\Omega}^2$.
In \cite{LangerRepinWolfmayr2016}, the corresponding optimality system is derived, which
is given in weak formulation as follows:
Given $y_d \in L^2(Q)$, find $y, p \in H^{1,\frac{1}{2}}_{0,per}(Q)$ such that
\begin{align}
 \label{problem:KKTSysSTVFAPost1}
  &\int_{Q} \Big( y\,z - \nu \nabla p \cdot \nabla z
      + \sigma \partial_t^{1/2} p \, \partial_t^{1/2} z^{\perp} \Big)\,d\boldsymbol{x}\,dt
      = \int_{Q} y_d\,z\,d\boldsymbol{x}\,dt, \qquad \forall z \in H^{1,\frac{1}{2}}_{0,per}(Q), 
      \\
   \label{problem:KKTSysSTVFAPost2}
  &\int_{Q} \Big( \nu \nabla y \cdot \nabla q
      + \sigma \partial_t^{1/2} y \, \partial_t^{1/2} q^{\perp}
      + \lambda^{-1} p\,q \Big)\,d\boldsymbol{x}\,dt = 0, \qquad \forall q \in H^{1,\frac{1}{2}}_{0,per}(Q).
\end{align}
The reduced optimality system \eqref{problem:KKTSysSTVFAPost1}--\eqref{problem:KKTSysSTVFAPost2} i.a.
results from using the condition $u = - \lambda^{-1} p$,
since no box constraints are imposed on the control function $u$ in this paper. 
This also leads to
the space of admissible controls being given by $H^{1,\frac{1}{2}}_{0,per}(Q)$.
\begin{remark}
\label{rem:fotd}
Since the optimality system is derived first, and then the discretization is performed later, we can say
that the \emph{first optimize, then discretize} approach is applied here (as discussed, e.g., earlier in
\cite{HinzeRoesch2012}).
\end{remark}

Let us define $V := H^1_0(\Omega)$ and $\mathbb{V} := V \times V$.
Expanding all functions into
Fourier series 
in \eqref{problem:KKTSysSTVFAPost1}--\eqref{problem:KKTSysSTVFAPost2} together with using  
the orthogonality of the cosine and sine 
functions yields the following problems
corresponding to the $k$th and 0th Fourier coefficients
and which are
decoupled due to the linearity of the optimal control problem:
Find $\boldsymbol{y}_k, \boldsymbol{p}_k \in \mathbb{V}$ 
so that
\begin{align}
 \label{equation:MultiAnsVFBlockk1}
 &\int_{\Omega} \big( \boldsymbol{y}_k \cdot \boldsymbol{z}_k
  - \nu \nabla \boldsymbol{p}_k \cdot \nabla \boldsymbol{z}_k
  + k \omega \, \sigma \boldsymbol{p}_k \cdot \boldsymbol{z}_k^{\perp} \big)\,d\boldsymbol{x}
  = \int_{\Omega} \boldsymbol{y_d}_k \cdot \boldsymbol{z}_k\,d\boldsymbol{x},
  \quad \forall \boldsymbol{z}_k \in \mathbb{V}, \\
  \label{equation:MultiAnsVFBlockk2}
 &\int_{\Omega} \big( \nu \nabla \boldsymbol{y}_k \cdot \nabla \boldsymbol{q}_k
  + k \omega \, \sigma \boldsymbol{y}_k \cdot \boldsymbol{q}_k^{\perp}
  + \lambda^{-1} \boldsymbol{p}_k \cdot \boldsymbol{q}_k \big)\,d\boldsymbol{x} = 0,
  \quad \forall \boldsymbol{q}_k \in \mathbb{V}.
\end{align}
and $y_0^c, p_0^c \in V$ so that
\begin{align}
 \label{equation:MultiAnsVFBlock01}
 &\int_{\Omega} \big( y_0^c \cdot z_0^c - \nu \nabla p_0^c \cdot \nabla z_0^c \big)\,d\boldsymbol{x}
 = \int_{\Omega} {y_d^c}_0 \cdot z_0^c\,d\boldsymbol{x}, 
 \quad \forall z_0^c \in V, \\
  \label{equation:MultiAnsVFBlock02}
 &\int_{\Omega} \big( \nu \nabla y_0^c \cdot \nabla q_0^c
 + \lambda^{-1} p_0^c \cdot q_0^c \big)\,d\boldsymbol{x} = 0,
 \quad \forall q_0^c \in V.
\end{align}
Both problems 
\eqref{equation:MultiAnsVFBlockk1}--\eqref{equation:MultiAnsVFBlockk2} 
and \eqref{equation:MultiAnsVFBlock01}--\eqref{equation:MultiAnsVFBlock02}
are uniquely solvable (see \cite{LangerWolfmayr:2013}).

\subsection{Majorant for 
cost functional \eqref{equation:minfunc1}}
\label{SSec:MajorantI}

Here, the results of \cite{LangerRepinWolfmayr2016} on upper bounds
for optimization problem I are summerized,
since they are needed later to derive the two-sided estimate for the 
cost functional \eqref{equation:minfunc1},
which deepens and extends the a posteriori error analysis
for optimization problem I.
Let $y=y(v)$ be the corresponding state to an arbitrary control $v$.
The following upper bound can be proved:
  \begin{align*}
   \mathcal{J}(y(v),v)
   \leq \mathcal{J}^\oplus(\alpha,\beta;\eta,\boldsymbol{\tau},v) \qquad \forall \, v \in L^2(Q),
  \end{align*}
for arbitrary  $\alpha, \beta > 0$, $\eta \in H^{1,1}_{0,per}(Q)$ and
\begin{align*}
 \boldsymbol{\tau} \in H(\text{div},Q) :=
 \{\boldsymbol{\tau} \in [L^2(Q)]^d 
 : \nabla \cdot \boldsymbol{\tau}(\cdot,t) \in L^2(\Omega) 
 \text{ for a.e. } t \in (0,T)\},
\end{align*}
where, for any $\boldsymbol{\tau} \in H(\text{div},Q)$, 
the identity
\begin{align*}
 \int_\Omega \nabla \cdot \boldsymbol{\tau} \, w \, d\boldsymbol{x} 
 = -  \int_\Omega \boldsymbol{\tau} \cdot \nabla w \, d\boldsymbol{x} \qquad \forall \, w \in V 
\end{align*}
is valid.
The guaranteed and fully computable majorant is given by
\begin{align}
  \label{definition:majorantCostFuncI}
  \begin{aligned}
   \mathcal{J}^\oplus(\alpha,\beta;\eta,\boldsymbol{\tau},v) :=
   &\, \frac{1+\alpha}{2} \|\eta - y_d\|^2
     + \gamma (
     \|\mathcal{R}_2(\eta,\boldsymbol{\tau})\|^2 
     + \frac{C_F^2}{ \beta}
      \|\mathcal{R}_1(\eta,\boldsymbol{\tau},v)\|^2 )
     + \frac{\lambda}{2} \|v\|^2,
  \end{aligned}
\end{align}
where $\underline{\mu_{1}} = \frac{1}{\sqrt{2}} \min\{\underline{\nu},\underline{\sigma}\}$
and we have set
$\gamma := \frac{(1+\alpha)(1+\beta) C_F^2 }{2 \alpha \underline{\mu_1}^2}$.
The parameters 
$\alpha, \beta > 0$ have been introduced in order
to obtain a quadratic functional by applying Young's inequality.
\begin{remark}
The arbitrary functions $\eta \in H^{1,1}_{0,per}(Q)$ and $v \in L^2(Q)$
can be taken later as the approximate solutions of the optimal control problem
\eqref{equation:minfunc1} subject to \eqref{equation:forwardpde:OCP1}--\eqref{equation:forwardpde:OCP3}
and $\boldsymbol{\tau} \in H(\emph{div},Q)$
represents the image of the flux $\nu \nabla \eta$. Note again that \emph{first optimize, then discretize}
is applied in this paper,
see also Remark \ref{rem:fotd}.
\end{remark}
For the derivation of \eqref{definition:majorantCostFuncI}, the
following estimate for the approximation error has been used:
 \begin{align}
  \label{inequality:aposteriorEstimateH11/2SeminormBVP}
  |y(v) - \eta|_{1,\frac{1}{2}} \leq \frac{1}{\underline{\mu_1}}
  \left(C_F \, \|\mathcal{R}_1(\eta,\boldsymbol{\tau},v)\| 
    + \|\mathcal{R}_2(\eta,\boldsymbol{\tau})\| \right),
 \end{align}
where
\begin{align}
\label{def:R1R2}
\mathcal{R}_1(\eta,\boldsymbol{\tau},v)
   := \sigma \partial_t \eta - \nabla \cdot \boldsymbol{\tau} - v \qquad \text{ and } \qquad
\mathcal{R}_2(\eta,\boldsymbol{\tau})
   := \boldsymbol{\tau} - \nu \nabla \eta.
\end{align}
The derivation of estimate \eqref{inequality:aposteriorEstimateH11/2SeminormBVP}
can be found in \cite{LangerRepinWolfmayr2015}.
The function $\mathcal{J}^\oplus(\alpha,\beta;\eta,\boldsymbol{\tau},v)$
is a sharp upper bound on $\mathcal{J}(y(v),v)$ for arbitrary but fixed $v$
as well as on
the optimal value $\mathcal{J}(y(u),u)$
 \begin{align}
  \label{definition:majorantCostFuncInfI}
   \inf_{\substack{\eta \in H^{1,1}_{0,per}(Q),\boldsymbol{\tau} \in H(\text{div},Q) \\
   v \in L^2(Q), \alpha, \beta > 0}}
   \mathcal{J}^\oplus(\alpha,\beta;\eta,\boldsymbol{\tau},v) = \mathcal{J}(y(u),u),
 \end{align}
since the infimum is attained for the optimal control $u$, its corresponding state $y(u)$
and its exact flux $\nu \nabla y(u)$, and for $\alpha$ going to zero.
Therefore, we have the estimate
  \begin{align}
  \label{estimate:majCostFuncI}
  \begin{aligned}
   \mathcal{J}(y(u),u)
   \leq \mathcal{J}^\oplus(\alpha,\beta;\eta,\boldsymbol{\tau},v) \quad
   \forall \, \eta \in H^{1,1}_{0,per}(Q), \, 
   \boldsymbol{\tau} \in H(\text{div},Q), \, 
   v \in L^2(Q), \, \alpha, \beta > 0.
  \end{aligned}
  \end{align}

\subsection{Optimization problem II}
\label{SSec:ModelProblemII}

In the second case, we want to minimize the following cost functional with respect to the unknown
state $y$ and control $u$:
\begin{align}
\label{equation:minfunc2}
 \tilde{\mathcal{J}}(y,u) :=  \frac{1}{2} \| \nabla y - \boldsymbol{g}_d \|^2 + \frac{\lambda}{2}  \| u \|^2
\end{align}
subject to the time-periodic boundary value problem \eqref{equation:forwardpde:OCP1}--\eqref{equation:forwardpde:OCP3},
where $\boldsymbol{g}_d \in [L^2(Q)]^d$ is the given desired gradient.
The optimality system can analogously be derived as for optimization problem I
using the Lagrange functional
\begin{align*}
   \tilde{\mathcal{L}}(y,u,p) = \tilde{\mathcal{J}}(y,u) 
   - \int_{Q} \left(  \sigma \, \partial_t y
 - \nabla \cdot ( \nu \nabla y) - u \right) p \, \dxdt.
\end{align*}
Optimality equations $\partial_u \tilde{\mathcal{L}}(y,u,p) = 0$ and $\partial_p \tilde{\mathcal{L}}(y,u,p) = 0$
are similar.
Equation $\partial_y \tilde{\mathcal{L}}(y,u,p) = 0$ is different. 
The optimality conditions are given in weak form as follows:
Given $\boldsymbol{g}_d \in [L^2(Q)]^d$, find $y, p \in
H^{1,\frac{1}{2}}_{0,per}(Q)$ such that
\begin{align}
 \label{problem:KKTSysSTVFAPostII1}
  &\int_{Q} \Big( 
  \nabla y \cdot \nabla 
  z - \nu \nabla p \cdot \nabla z
      + \sigma \partial_t^{1/2} p \, \partial_t^{1/2} z^{\perp} \Big)\,d\boldsymbol{x}\,dt
      = \int_{Q} 
      \boldsymbol{g}_d \cdot 
      \nabla z\,d\boldsymbol{x}\,dt, \quad \forall z \in H^{1,\frac{1}{2}}_{0,per}(Q), \\
    \label{problem:KKTSysSTVFAPostII2}
  &\int_{Q} \Big( \nu \nabla y \cdot \nabla q
      + \sigma \partial_t^{1/2} y \, \partial_t^{1/2} q^{\perp}
      + \lambda^{-1} p\,q \Big)\,d\boldsymbol{x}\,dt = 0, \quad \forall q \in H^{1,\frac{1}{2}}_{0,per}(Q).
\end{align}
The optimality systems corresponding to every mode $k$ are analogously
derived as for optimization problem I 
(similar to \eqref{equation:MultiAnsVFBlockk1}--\eqref{equation:MultiAnsVFBlockk2}
and \eqref{equation:MultiAnsVFBlock01}--\eqref{equation:MultiAnsVFBlock02}).
In Section \ref{Sec:TwoSidedBoundsII}, we will derive new two-sided bounds for
optimization problem II. 

\section{Guaranteed lower bounds leading to 
two-sided bounds for optimization problem I}
\label{Sec:TwoSidedBoundsI}

In this work, we complement the  
guaranteed upper bounds for the discretization error
in state and control of minimizing cost functional $\mathcal{J}$
defined in
\eqref{equation:minfunc1} subject to \eqref{equation:forwardpde:OCP1}--\eqref{equation:forwardpde:OCP3}.
This is done
by obtaining fully computable lower bounds for 
$\mathcal{J}$ following the technique from
\cite{Wolfmayr2016} (derived for elliptic problems)
leading to two-sided bounds for the cost functional \eqref{equation:minfunc1}.

\subsection{Minorant for cost functional \eqref{equation:minfunc1}}
\label{SSec:MinorantI}

Let $y=y(u)$ be the optimal 
state corresponding to the optimal control function 
$u \in L^2(Q)$, which is connected with the adjoint state $p = p(u)$
by the identity $u = - \lambda^{-1} p(u)$.
Then $y=y(u)$ is the solution of the variational formulation
\begin{align}
  \label{equation:forwardpde:VF}
 \int_{Q} \Big( \nu \nabla y \cdot \nabla q
      + \sigma \partial_t^{1/2} y \, \partial_t^{1/2} q^{\perp}
      + \lambda^{-1} p\,q \Big) \dxdt = 0 \qquad \forall \, q \in H^{1,\frac{1}{2}}_{0,per}(Q)
\end{align}
of the boundary value problem \eqref{equation:forwardpde:OCP1}--\eqref{equation:forwardpde:OCP3} (see also
\eqref{problem:KKTSysSTVFAPost1}--\eqref{problem:KKTSysSTVFAPost2}).
For any arbitray function $\eta \in H^{1,1}_{0,per}(Q)$, 
one can obtain that
\begin{align*}
\mathcal{J}(y(u),u) =  \frac{1}{2} \| y - \eta \|^2
+ \int_{Q}
\left(y-\eta\right)\left(\eta-y_d\right) \dxdt
+ \frac{1}{2} \| \eta - y_d \|^2 + \frac{\lambda}{2}  \| u \|^2.
\end{align*}
Since $\frac{1}{2} \|y-\eta\|^2 \geq 0$ and
using the identity $u = - \lambda^{-1} p(u)$, $\mathcal{J}$ can be estimated from below by
\begin{align}
\label{equation:costfunctionalExpandedEstBelow}
\mathcal{J}(y(u),u) = 
  \mathcal{J}(y(u),- \lambda^{-1} p(u))
  \geq \frac{1}{2} \|\eta-y_d\|^2
 + \frac{1}{2 \lambda} \|p\|^2 + \int_{Q}
 \left(y-\eta\right)\left(\eta-y_d\right) \dxdt.
\end{align}
For $\eta \in H^{1,1}_{0,per}(Q)$,
let $p_\eta,  \tilde{p}_\eta \in H^{1,\frac{1}{2}}_{0,per}(Q)$ be the 
solutions to the equations
\begin{align}
 \label{equation:NOCforPeta1}
   &\int_{Q} \Big(  \nu \nabla p_\eta \cdot \nabla z
      - \sigma \partial_t^{1/2} p_\eta \, \partial_t^{1/2} z^{\perp} \Big)\dxdt
      = \int_{Q} (\eta - y_d) z \, \dxdt, \qquad \forall z \in H^{1,\frac{1}{2}}_{0,per}(Q), \\
 \label{equation:NOCforPeta2}
    &\int_{Q} \Big( \nu \nabla \eta \cdot \nabla q
      + \sigma \partial_t^{1/2} \eta \, \partial_t^{1/2} q^{\perp}
      + \lambda^{-1} \tilde{p}_\eta \,q \Big)\,d\boldsymbol{x}\,dt = 0, 
      \qquad \forall q \in H^{1,\frac{1}{2}}_{0,per}(Q).
\end{align}
\begin{remark}
\label{Rem1}
Note that we assumed that $\eta \in H^{1,1}_{0,per}(Q)$
according to the derivation of the majorant, but so far the assumption
$\eta \in H^{1,\frac{1}{2}}_{0,per}(Q)$ would be enough.
\end{remark}
Adding and subtracting $p_\eta$ in \eqref{equation:costfunctionalExpandedEstBelow}
together with $\frac{1}{2 \lambda} \|p - p_\eta\|^2 \geq 0$
yields the estimate
\begin{align*}
\mathcal{J}(y(u),u) 
  &\geq \frac{1}{2} \|\eta-y_d\|^2
 + \frac{1}{2 \lambda} \|p_\eta\|^2 
 + \intQT \left(y-\eta\right)\left(\eta-y_d\right) \dxdt
 + \intQT \lambda^{-1} \left(p - p_\eta\right) p_\eta \, \dxdt .
\end{align*}
Using equation \eqref{equation:NOCforPeta1} and
identity \eqref{def:identityH01/2per} leads to the estimate 
\begin{align}
\label{est:minfunc1long}
\left \lbrace 
\begin{aligned}
 \mathcal{J}(y(u),u) \geq 
 &\, \frac{1}{2} \|\eta-y_d\|^2
 +  \frac{1}{2 \lambda} \|p_\eta\|^2 
 + \intQT \lambda^{-1} \left(p - p_\eta\right) p_\eta \, \dxdt \\
  &+\int_{Q} \Big(  \nu \nabla p_\eta \cdot \nabla \left(y-\eta\right) 
    - \sigma \partial_t^{1/2} p_\eta \, \partial_t^{1/2} \left(y-\eta\right)^{\perp} \Big)\dxdt \\
 = &\, \frac{1}{2} \|\eta-y_d\|^2
 +  \frac{1}{2 \lambda} \|p_\eta\|^2 
 + \intQT \lambda^{-1} \left(p - p_\eta\right) p_\eta \, \dxdt \\
&+ \intQT \Big( \left(\nu \nabla y - \nu \nabla \eta\right) \cdot \nabla {p_\eta} 
 + \left(\sigma \partial_t^{1/2} y- \sigma \partial_t^{1/2} \eta\right) \partial_t^{1/2} p_\eta^{\perp} \Big) \dxdt \\ 
 = &\, \frac{1}{2} \|\eta-y_d\|^2
 +  \frac{1}{2 \lambda} \|p_\eta\|^2 
+ \intQT \Big( \nu \nabla y \cdot \nabla {p_\eta} 
 + \sigma \partial_t^{1/2} y \, \partial_t^{1/2} p_\eta^{\perp}
 + \lambda^{-1} p \, p_\eta \Big) \dxdt \\
&- \intQT \Big(\nu \nabla \eta \cdot \nabla {p_\eta} 
 + \sigma \partial_t^{1/2} \eta \, \partial_t^{1/2} p_\eta^{\perp}
 + \lambda^{-1} p_\eta \, p_\eta \Big) \dxdt.
\end{aligned}
\right.
\end{align}
By applying
equations \eqref{equation:forwardpde:VF} and \eqref{equation:NOCforPeta2},
it follows that
\begin{align*}
 \mathcal{J}(y(u),u) 
\geq &\, \frac{1}{2} \|\eta-y_d\|^2
 +  \frac{1}{2 \lambda} \|p_\eta\|^2
 + \intQT \lambda^{-1} (\tilde{p}_\eta - p_\eta) p_\eta \, \dxdt. 
\end{align*}
We introduce now
the arbitrary function $\zeta \in H^{1,1}_{0,per}(Q)$.
Note that at the moment
$\zeta \in H^{1,\frac{1}{2}}_{0,per}(Q)$ would be enough,
but the higher regularity in time will be needed in another step.
This goes along with the higher regularity assumption on $\eta$ (see
Remark \ref{Rem1}).
Since $\frac{1}{2 \lambda} \|p_\eta-\zeta\|^2 \geq 0$, we have that
\begin{align*}
 \mathcal{J}(y(u),u) \geq 
&\, \frac{1}{2} \|\eta-y_d\|^2
 +  \frac{1}{2 \lambda} \|\zeta\|^2
 + \intQT \lambda^{-1} \left(p_\eta \zeta -\zeta^2 + \tilde{p}_\eta p_\eta - p_\eta^2 \right) \, \dxdt.
\end{align*}
Now we add and subtract $\lambda^{-1} \tilde{p}_\eta \zeta$ in the last integral
as well as use equation \eqref{equation:NOCforPeta2} again.
Moreover, we exploit the fact that we have 
assumed that $\eta \in H^{1,1}_{0,per}(Q)$, hence, 
we can apply
the identities \eqref{equation:H11/2identities}.
Altogether these steps
yield the estimate
\begin{align}
\label{preest}
\begin{aligned}
 \mathcal{J}(y(u),u) \geq 
 &\, \frac{1}{2} \|\eta-y_d\|^2
 +  \frac{1}{2 \lambda} \|\zeta\|^2
 - \intQT \left(\nu \nabla \eta \cdot \nabla \zeta
      + \sigma \partial_t \eta \, \zeta
      + \lambda^{-1} \zeta^2 \right) \, \dxdt \\
 &+ \intQT \lambda^{-1} (\zeta - p_\eta) (p_\eta-\tilde{p}_\eta) \, \dxdt.
\end{aligned}
\end{align}
In the following, we need to estimate the last integral of this expression in
order to formulate a computable lower bound for the cost functional.
For that let us first prove a computable
upper bound for the error in the adjoint state, which is presented in the following theorem.
Note that here we will need the higher regularity assumption (in time) on $\zeta$.
\begin{theorem}
\label{theorem:estimate:adjointState}
 Let $y_d \in L^2(Q)$ be given. Let 
 $p_\eta \in H^{1,\frac{1}{2}}_{0,per}(Q)$ solve
  \eqref{equation:NOCforPeta1} for an arbitrary 
 $\eta \in H^{1,1}_{0,per}(Q)$.
 The following estimate holds:
 \begin{align}
  \label{estimate:adjointState}
  \|\nabla ({p_\eta}- \zeta)\| \leq \frac{1}{\underline{\mu_{1}}}
  \left(C_F \|\mathcal{R}_3(\zeta,\boldsymbol{\rho},\eta)\| + \|\mathcal{R}_4(\zeta,\boldsymbol{\rho})\| \right)
 \end{align}
 for any $\zeta \in H^{1,1}_{0,per}(Q)$,
 where $\underline{\mu_{1}} = \frac{1}{\sqrt{2}} \min\{\underline{\nu},\underline{\sigma}\}$,
 $\mathcal{R}_3(\zeta,\boldsymbol{\rho},\eta) 
 = \eta - y_d + \nabla \cdot \boldsymbol{\rho} + \sigma \partial_t \zeta$ and
 $\mathcal{R}_4(\zeta,\boldsymbol{\rho}) = \boldsymbol{\rho} - \nu \nabla \zeta$
 with $\boldsymbol{\rho} \in H(\emph{div},Q)$
 and $C_F > 0$. 
 \begin{proof}
 Let us consider the adjoint equation \eqref{equation:NOCforPeta1}.
 Adding and subtracting  $\zeta \in H^{1,1}_{0,per}(Q)$ in the left-hand side of the equation leads to
 \begin{align}
 \label{auxeq}
 \begin{aligned}
   \int_{Q} &\Big(  \nu \nabla (p_\eta - \zeta) \cdot \nabla z
      - \sigma \partial_t^{1/2} (p_\eta - \zeta) \, \partial_t^{1/2} z^{\perp} \Big)\dxdt
      = \int_{Q} (\eta - y_d) z \, \dxdt \\
         &- \int_{Q} \nu \nabla \zeta \cdot \nabla z \, \dxdt
      + \int_{Q} \sigma \partial_t^{1/2} \zeta \, \partial_t^{1/2} z^{\perp} \dxdt.
\end{aligned}
\end{align}
Next we introduce the auxiliary variable $\boldsymbol{\rho} \in H(\text{div},Q)$.
Together with using that $\zeta \in H^{1,1}_{0,per}(Q)$ as well as applying
Cauchy--Schwarz' and Friedrichs' inequalities, the
following estimate for the right-hand side of \eqref{auxeq} can be obtained:
 \begin{align*}
  &\sup_{0 \not=z \in H^{1,\frac{1}{2}}_{0,per}(Q)} \frac{\int_{Q} (\eta - y_d + \nabla \cdot \boldsymbol{\rho}
  + \sigma \partial_t \zeta) z \, \dxdt 
         + \int_{Q}(\boldsymbol{\rho} - \nu \nabla \zeta) \cdot \nabla z \, \dxdt}{|z|_{1,\frac{1}{2}}
         } \\
   &\leq \sup_{0 \not=z \in H^{1,\frac{1}{2}}_{0,per}(Q)} \frac{\|\eta - y_d + \nabla \cdot \boldsymbol{\rho}
  + \sigma \partial_t \zeta \| \| z \|
         + \|\boldsymbol{\rho} - \nu \nabla \zeta \| \| \nabla z \|}{|z|_{1,\frac{1}{2}}
         } \\
   &\leq \sup_{0 \not=z \in H^{1,\frac{1}{2}}_{0,per}(Q)} \frac{C_F \|\eta - y_d + \nabla \cdot \boldsymbol{\rho}
  + \sigma \partial_t \zeta \| \| \nabla z \|
         + \|\boldsymbol{\rho} - \nu \nabla \zeta \| \| \nabla z \|}{\| \nabla z\|} \\
   &\leq C_F \|\eta - y_d + \nabla \cdot \boldsymbol{\rho}
  + \sigma \partial_t \zeta \|
         + \|\boldsymbol{\rho} - \nu \nabla \zeta \|.
\end{align*}
Using the boundedness of the coefficients $\sigma$ and $\nu$,
the orthogonality relations
\eqref{equation:orthorelation2} and applying that $\|z^\perp\| = \|z\|$,
we can prove the  
estimate from below for the left-hand side of 
\eqref{auxeq}, which is
 \begin{align*}
  \sup_{0 \not=z \in H^{1,\frac{1}{2}}_{0,per}(Q)} \frac{\int_{Q} \Big(  \nu \nabla (p_\eta - \zeta) \cdot \nabla z
      - \sigma \partial_t^{1/2} (p_\eta - \zeta) \, \partial_t^{1/2} z^{\perp} \Big)\dxdt}{|z|_{1,\frac{1}{2}} }.
\end{align*}
First, we estimate the supremum from below with choosing $z = (p_\eta - \zeta) + (p_\eta - \zeta)^\perp$,
for which
 \begin{align*}
 |z|_{1,\frac{1}{2}} = |(p_\eta - \zeta) + (p_\eta - \zeta)^\perp|_{1,\frac{1}{2}} = \sqrt{2} \, |p_\eta - \zeta|_{1,\frac{1}{2}},
\end{align*}
using the orthogonality relations \eqref{equation:orthorelation2}.
Next, applying the second equation in \eqref{equation:orthorelation2} gives
 \begin{align*}
 \langle \nu \nabla (p_\eta - \zeta), \nabla z \rangle
  &= \langle \nu \nabla (p_\eta - \zeta), \nabla ((p_\eta - \zeta) + (p_\eta - \zeta)^\perp) \rangle \\
  &= \langle \nu \nabla (p_\eta - \zeta), \nabla (p_\eta - \zeta)\rangle + \langle \nu \nabla (p_\eta - \zeta), \nabla (p_\eta - \zeta)^\perp \rangle \\
  &= \langle \nu \nabla (p_\eta - \zeta), \nabla (p_\eta - \zeta)\rangle,
\end{align*}
and applying the first equation in \eqref{equation:orthorelation2} as well as using the identity $((p_\eta - \zeta)^\perp)^{\perp} = - (p_\eta - \zeta)$ gives
 \begin{align*}
 \langle \sigma \partial_t^{1/2} (p_\eta - \zeta), \partial_t^{1/2} z^{\perp} \rangle 
 &= \langle \sigma \partial_t^{1/2} (p_\eta - \zeta), \partial_t^{1/2} ((p_\eta - \zeta) + (p_\eta - \zeta)^\perp)^{\perp} \rangle \\
  &=  \langle \sigma \partial_t^{1/2} (p_\eta - \zeta), \partial_t^{1/2} (p_\eta - \zeta)^{\perp} \rangle
  +  \langle \sigma \partial_t^{1/2} (p_\eta - \zeta), \partial_t^{1/2} ((p_\eta - \zeta)^\perp)^{\perp} \rangle \\
  &= - \langle \sigma \partial_t^{1/2} (p_\eta - \zeta), \partial_t^{1/2} (p_\eta - \zeta) \rangle
\end{align*}
leading to the estimate
 \begin{align*}
  &\sup_{0 \not=z \in H^{1,\frac{1}{2}}_{0,per}(Q)} \frac{\int_{Q} \Big(  \nu \nabla (p_\eta - \zeta) \cdot \nabla z
      - \sigma \partial_t^{1/2} (p_\eta - \zeta) \, \partial_t^{1/2} z^{\perp} \Big)\dxdt}{|z|_{1,\frac{1}{2}}
      } \\
   &\geq \frac{\langle \nu \nabla (p_\eta - \zeta), \nabla (p_\eta - \zeta)\rangle
      + \langle \sigma \partial_t^{1/2} (p_\eta - \zeta), \partial_t^{1/2} (p_\eta - \zeta) \rangle}{|(p_\eta - \zeta) + (p_\eta - \zeta)^\perp|_{1,\frac{1}{2}}
      }\\
    &= \frac{\langle \nu \nabla (p_\eta - \zeta), \nabla (p_\eta - \zeta)\rangle
      + \langle \sigma \partial_t^{1/2} (p_\eta - \zeta), \partial_t^{1/2} (p_\eta - \zeta) \rangle}{\sqrt{2} \, |p_\eta - \zeta|_{1,\frac{1}{2}} }\\
    &\geq \frac{\underline{\nu} \| \nabla (p_\eta - \zeta)\|^2
      +\underline{\sigma} \| \partial_t^{1/2} (p_\eta - \zeta)\|^2}{\sqrt{2} \, |p_\eta - \zeta|_{1,\frac{1}{2}} }\\
   &\geq \underline{\mu_{1}} \, |p_\eta - \zeta|_{1,\frac{1}{2}},
\end{align*}
where $\underline{\mu_{1}} = \frac{1}{\sqrt{2}} \min\{\underline{\nu},\underline{\sigma}\}$.
Combining now both estimates together with
$\|\nabla (p_\eta - \zeta)\| \leq |p_\eta - \zeta|_{1,\frac{1}{2}}$
we finally derive 
\eqref{estimate:adjointState}.
 \end{proof}
\end{theorem}
Now we have all the tools in order to estimate the last term of \eqref{preest}.
We obtain as follows
\begin{align*}
& \intQT \lambda^{-1} (\zeta - p_\eta) (p_\eta - \zeta + \zeta -\tilde{p}_\eta) \, \dxdt 
 = \intQT (\lambda^{-1} (\zeta - p_\eta) (p_\eta - \zeta)
 + \lambda^{-1} (\zeta - p_\eta) (\zeta -\tilde{p}_\eta)) \, \dxdt \\
 &= - \lambda^{-1} \|\zeta - p_\eta\|^2
 + \intQT ( \lambda^{-1} (\zeta - p_\eta) (\zeta -\tilde{p}_\eta)) \, \dxdt \\
 &= - \lambda^{-1} \|\zeta - p_\eta\|^2
 + \int_{Q} \Big( \nu \nabla \eta \cdot \nabla (\zeta - p_\eta)
      + \sigma \partial_t^{1/2} \eta \, \partial_t^{1/2} (\zeta - p_\eta)^{\perp} \Big)\,d\boldsymbol{x}\,dt 
      + \int_{Q}  \lambda^{-1} \zeta (\zeta - p_\eta) \,d\boldsymbol{x}\,dt \\
  &= - \lambda^{-1} \|\zeta - p_\eta\|^2
 + \int_{Q} \left(\sigma \partial_t \eta
      - \nabla \cdot \boldsymbol{\tau}
      + \lambda^{-1} \zeta \right) (\zeta - p_\eta) \,d\boldsymbol{x}\,dt 
 + \int_{Q} \left( \nu \nabla \eta
      - \boldsymbol{\tau}\right) \cdot \nabla (\zeta - p_\eta) \,d\boldsymbol{x}\,dt \\
  &\geq 
 - \lambda^{-1} C_F^2 \|\nabla(\zeta - p_\eta)\|^2
 - (C_F \|\mathcal{R}_1(\eta,\boldsymbol{\tau},-\lambda^{-1} \zeta)\|
 + \|\mathcal{R}_2(\eta,\boldsymbol{\tau})\|) \| \nabla (\zeta - p_\eta) \| 
\end{align*}
leading to
\begin{align}
\label{preest2}
\begin{aligned}
& \intQT \lambda^{-1} (\zeta - p_\eta) (p_\eta - \zeta + \zeta -\tilde{p}_\eta) \, \dxdt 
   \geq  - \frac{C_F^2}{\underline{\mu_{1}}^2 \lambda}
  \left(C_F \|\mathcal{R}_3(\zeta,\boldsymbol{\rho},\eta)\|
  + \|\mathcal{R}_4(\zeta,\boldsymbol{\rho})\| \right)^2 \\
 &\qquad - \frac{1}{\underline{\mu_{1}}} 
 (C_F \|\mathcal{R}_1(\eta,\boldsymbol{\tau},-\lambda^{-1} \zeta)\| 
      + \| \mathcal{R}_2(\eta,\boldsymbol{\tau})\|) 
  \left(C_F \|\mathcal{R}_3(\zeta,\boldsymbol{\rho},\eta)  \| 
  + \|\mathcal{R}_4(\zeta,\boldsymbol{\rho}) \| \right),
 \end{aligned}
\end{align}
where  $\boldsymbol{\tau}, \boldsymbol{\rho} \in H(\text{div},Q)$
and we have used equation \eqref{equation:NOCforPeta2},
relations \eqref{equation:orthorelation1}--\eqref{equation:orthorelation2},
Cauchy--Schwarz' and 
Friedrichs' inequalities, estimate \eqref{estimate:adjointState}
and that $\eta \in H^{1,1}_{0,per}(Q)$.

Finally, we obtain the following estimate from below:
\begin{align}
\label{inequality:finalminorant1}
 \mathcal{J}(y(u),u) \geq 
 \mathcal{J}^\ominus(\eta,{\zeta},\boldsymbol{\tau},\boldsymbol{\rho}) \qquad
 \forall \, \eta, \zeta \in H^{1,1}_{0,per}(Q), 
  \boldsymbol{\tau},\boldsymbol{\rho} 
  \in H(\text{div},Q),
\end{align}
where the (fully computable) minorant is given by
\begin{align}
 \label{definition:minorant}
\begin{aligned}
 \mathcal{J}^\ominus(\eta,{\zeta}&,\boldsymbol{\tau},\boldsymbol{\rho})
 =  \frac{1}{2} \|\eta-y_d\|^2
 +  \frac{1}{2 \lambda} \|\zeta\|^2
 - \int_{Q} 
 \left(\nu \nabla \eta \cdot \nabla \zeta
      + \sigma \partial_t \eta \, \zeta
      + \lambda^{-1} \zeta^2 \right) \, \dxdt \\
&- \frac{C_F^2}{\underline{\mu_{1}}^2 \lambda}
  \left(C_F \|\mathcal{R}_3(\zeta,\boldsymbol{\rho},\eta)\| 
  + \|\mathcal{R}_4(\zeta,\boldsymbol{\rho})\| \right)^2 \\
 & - \frac{1}{\underline{\mu_{1}}} (C_F \|\mathcal{R}_1(\eta,\boldsymbol{\tau},-\lambda^{-1} \zeta)\| 
      + \|\mathcal{R}_2(\eta,\boldsymbol{\tau})\|) 
  \left(C_F \|\mathcal{R}_3(\zeta,\boldsymbol{\rho},\eta)\| 
  + \|\mathcal{R}_4(\zeta,\boldsymbol{\rho})\| \right).
\end{aligned}
\end{align}
\begin{theorem} 
\label{theorem2}
 The supremum 
 of the minorant $\mathcal{J}^\ominus$ as
 defined in \eqref{definition:minorant} is attained for  
 the minimum  
 of the cost functional 
 \eqref{equation:minfunc1} subject to
 \eqref{equation:forwardpde:OCP1}--\eqref{equation:forwardpde:OCP3}, which is equivalent
 to the optimal value of 
 the optimality system \eqref{problem:KKTSysSTVFAPost1}--\eqref{problem:KKTSysSTVFAPost2}
 as follows
\begin{align}
 \label{definition:minorantSup}
  \sup_{\substack{\eta, \zeta \in H^{1,1}_{0,per}(Q), 
  \boldsymbol{\tau},\boldsymbol{\rho} 
  \in H(\emph{div},Q)}}
  \mathcal{J}^\ominus(\eta,{\zeta},\boldsymbol{\tau},\boldsymbol{\rho})
  = \mathcal{J}(y(u),u).
\end{align}
\begin{proof}
The estimate is sharp
for the exact control $u$, 
$\eta = y(u)$, $\zeta= p(u)$, $\boldsymbol{\tau} = \nu \nabla y(u)$
and $\boldsymbol{\rho} = \nu \nabla p(u)$. Hence,
 $\mathcal{J}^\ominus(y(u),p(u),\nu \nabla y(u)),\nu \nabla p(u))
 = \frac{1}{2} \|y-y_d\|^2
 + \frac{1}{2 \lambda} \|p\|^2
 = \frac{1}{2} \|y-y_d\|^2
 + \frac{\lambda}{2} \|u\|^2 = \mathcal{J}(y(u),u)$.
\end{proof}
\end{theorem}
\begin{remark}
\label{remark:convergence}
It has been shown in \cite{GaevskayaHoppeRepin:2006} that
if we choose finite dimensional subspaces that are limit dense
in the spaces of the exact solution $y(u)$, of its adjoint $p(u)$ and of their fluxes, $\nu \nabla y(u)$ and $\nu \nabla p(u)$,
which are
$H^{1,1}_{0,per}(Q)$ and $H(\emph{div},Q)$,
and we choose sequences of the functions $(\eta,{\zeta},\boldsymbol{\tau},\boldsymbol{\rho})$ in these
finite dimensional subspaces,
so for instance $(\eta_i,{\zeta}_i,\boldsymbol{\tau}_i,\boldsymbol{\rho}_i)$, then
they converge (let $i \rightarrow \infty$) to the exact solution, its adjoint and their fluxes.
The corresponding majorants $\mathcal{J}^\oplus_i$ and minorants $\mathcal{J}^\ominus_i$
converge to the exact value of the cost functional $\mathcal{J}$.
As shown in \cite{GaevskayaHoppeRepin:2006}, majorants
can be used in order to produce sequences of state and controls with values of the cost functional being
as close to the exact cost functional value as one needs it.
\end{remark}

\subsection{Guaranteed upper bounds for the discretization errors of the control 
and the state}
\label{SSec:MajControlStateI}

Here, we present
a posteriori error estimates for control and state
measured by the norm 
$ |||u-v|||^2 := \frac{1}{2} \|y(u) - y(v)\|^2 + \frac{\lambda}{2} \|u-v\|^2$.
The next theorem was proved for the elliptic case (together with control constraints) in 
\cite{Wolfmayr2016}. 
The norm $|||\cdot|||$ 
can be represented in terms of the state and the adjoint state (instead of the control), since
there are no control constraints imposed. Hence, 
$u = - \lambda^{-1} p(u)$, $v = - \lambda^{-1} p(v)$, 
and
 $|||u-v|||^2 = \frac{1}{2} \|y(u) - y(v)\|^2 + \frac{1}{2\lambda} \|p(u) - p(v)\|^2$.
\begin{theorem}
 \label{theorem:DiffCombNorm}
 We obtain the following identity for an arbitrary control $v \in L^2(Q)$:
 \begin{align}
  \label{equation:DiffCombNorm}
  |||u-v|||^2 = \mathcal{J}(y(v),v) - \mathcal{J}(y(u),u).
 \end{align}
\begin{proof}
Together with $u = - \lambda^{-1} p(u)$ and $v = - \lambda^{-1} p(v)$ the difference can be computed as
 \begin{align*}
  \mathcal{J}(y(v)&,v) - \mathcal{J}(y(u),u)
  = \frac{1}{2} \|y(v)-y_d\|^2 - \frac{1}{2} \|y(u)-y_d\|^2
    + \frac{\lambda}{2} \|v\|^2 - \frac{\lambda}{2} \|u\|^2 \\
  = &\,\frac{1}{2} \int_{Q} (y(v)-y(u)+2 y(u)-2 y_d)(y(v)-y(u)) \, \dxdt
    + \frac{\lambda}{2} \int_{Q} (v-u+2 u)(v-u) \, \dxdt \\
  = &\,\frac{1}{2} \|y(u) - y(v)\|^2 + \int_{Q} (y(u)-y_d)(y(v)-y(u)) \, \dxdt \\
    &+ \frac{\lambda}{2} \|u-v\|^2 + \lambda^{-1} \int_{Q} p(u)(p(v)-p(u)) \, \dxdt.
 \end{align*}
 Since the adjoint states $p(u),p(v) \in H^{1,1}_{0,per}(Q)$ fulfill
 \eqref{problem:KKTSysSTVFAPost1}--\eqref{problem:KKTSysSTVFAPost2} for the corresponding states
 $y(u)$, $y(v) \in H^{1,1}_{0,per}(Q)$, we obtain 
 \begin{align*}
  \mathcal{J}(&y(v),v) - \mathcal{J}(y(u),u) 
  =
  \frac{1}{2} \|y(u) - y(v)\|^2 + \frac{\lambda}{2} \|u-v\|^2 
    + \lambda^{-1} \int_{Q} p(u)(p(v)-p(u)) \, \dxdt \\
  & + \int_{Q} (\nu \nabla p(u) (\nabla y(v)- \nabla y(u)) 
     - \sigma \partial_t^{1/2} p(u) \, \partial_t^{1/2} (y(v) - y(u))^{\perp} )\, \dxdt 
  = \frac{1}{2} \|y(u) - y(v)\|^2 \\
  &+ \frac{\lambda}{2} \|u-v\|^2 
     + \lambda^{-1} \int_{Q} p(u)p(v) \, \dxdt 
    + \int_{Q} (\nu \nabla y(v) \cdot \nabla p(u) 
     + \sigma \partial_t^{1/2} y(v) \, \partial_t^{1/2} p(u)^{\perp}) 
     \, \dxdt \\
   &=  \frac{1}{2} \|y(u) - y(v)\|^2 
  + \frac{\lambda}{2} \|u-v\|^2.
 \end{align*}
 This proves now 
 the equality relation \eqref{equation:DiffCombNorm}
 by applying the equations \eqref{problem:KKTSysSTVFAPost1}--\eqref{problem:KKTSysSTVFAPost2} 
 for $(u,y(u),p(u))$ as well as $(v,y(v),p(v))$.
\end{proof}
\end{theorem}
Using the result of Theorem \ref{theorem:DiffCombNorm}, we can derive the majorant for
the discretization errors of control and state 
in the norm $|||\cdot|||$.
\begin{theorem}
\label{theorem:majorantCombNorm}
The functional
\begin{align*}
 \mathcal{M}^\oplus(\alpha,\beta&;\eta,\zeta,\boldsymbol{\tau},\boldsymbol{\rho},v) =
   \frac{\alpha}{2} \|\eta - y_d\|^2
     + \gamma (
     \|\mathcal{R}_2(\eta,\boldsymbol{\tau})\|^2 
     + \frac{C_F^2}{\beta}
      \|\mathcal{R}_1(\eta,\boldsymbol{\tau},v)\|^2 )
     + \frac{\lambda}{2} \|v\|^2 - \frac{1}{2 \lambda} \|\zeta\|^2 \\
   & + \intQT \left(\nu \nabla \eta \cdot \nabla \zeta
      + \sigma \partial_t \eta \, \zeta
      + \lambda^{-1} \zeta^2 \right) \, \dxdt + \frac{C_F^2}{\underline{\mu_1}^2 \lambda}
  \left(C_F \|\mathcal{R}_3(\zeta,\boldsymbol{\rho},\eta)\| 
  + \|\mathcal{R}_4(\zeta,\boldsymbol{\rho})\| \right)^2 \\
 &+ \frac{1}{\underline{\mu_1}} (C_F 
 \|\mathcal{R}_1(\eta,\boldsymbol{\tau},-\lambda^{-1} \zeta)\| 
      + \|\mathcal{R}_2(\eta,\boldsymbol{\tau})\|) 
  \left(C_F \|\mathcal{R}_3(\zeta,\boldsymbol{\rho},\eta)\| 
  + \|\mathcal{R}_4(\zeta,\boldsymbol{\rho})\| \right)
\end{align*}
for an arbitrary $v \in L^2(Q)$,
$\eta, \zeta \in H^{1,1}_{0,per}(Q)$,
$\boldsymbol{\tau}, \boldsymbol{\rho}
 \in H(\emph{div},Q)$
and 
$\alpha,\beta > 0$, and where
$\underline{\mu_{1}} = \frac{1}{\sqrt{2}} \min\{\underline{\nu},\underline{\sigma}\}$,
is a majorant for the discretization error
\begin{align}
 \label{definition:majorantCombinedNorm}
  |||u-v|||^2 \leq \mathcal{M}^\oplus(\alpha,\beta;\eta,\zeta,\boldsymbol{\tau},\boldsymbol{\rho},v)
  = \mathcal{J}^\oplus(\alpha,\beta;\eta,\boldsymbol{\tau},v)
  - \mathcal{J}^\ominus(\eta,{\zeta},\boldsymbol{\tau},\boldsymbol{\rho}).
\end{align}
\begin{proof}
  Applying \eqref{equation:DiffCombNorm} together with  
  \eqref{estimate:majCostFuncI} and \eqref{inequality:finalminorant1}
 yields estimate \eqref{definition:majorantCombinedNorm}.
\end{proof}
\end{theorem}
\begin{proposition}
\label{proposition:overallmajorant}
 The infimum of the majorant
 $\mathcal{M}^\oplus(\alpha,\beta;\eta,\zeta,\boldsymbol{\tau},\boldsymbol{\rho},v)$ 
 in \eqref{definition:majorantCombinedNorm} is attained 
 for the minimum  
 of the optimization problem I,
 which is equivalent to the solution of
 the optimality system \eqref{problem:KKTSysSTVFAPost1}--\eqref{problem:KKTSysSTVFAPost2}
 $(v=u$, $\eta = y(u)$, ${\zeta}= p(u) = - \lambda u$,
 $\boldsymbol{\tau} = \nu \nabla y(u)$,
 $\boldsymbol{\rho} = \nu \nabla p(u))$
 as follows
\begin{align*}
    \inf_{\substack{\eta,\zeta \in H^{1,1}_{0,per}(Q), 
    \boldsymbol{\tau}, \boldsymbol{\rho} \in H(\emph{div},Q), \\
     v \in L^2(Q), \alpha,\beta > 0}}
  \mathcal{M}^\oplus(\alpha,\beta;\eta,\zeta,\boldsymbol{\tau},\boldsymbol{\rho},v) = 0.
\end{align*}
\begin{proof}
 The majorant
$ \mathcal{M}^\oplus(\alpha,\beta;y(u),p(u),\nu \nabla y(u),\nu \nabla p(u),u)
 = \, \frac{\alpha}{2}\|y(u) - y_d\|^2$ 
 equals zero if $\alpha$ goes to zero.
\end{proof}
\end{proposition}

The majorant $\mathcal{M}^{\oplus}$ is a sharp, guaranteed and fully
 computable upper bound 
 for the control-state 
 error in $|||\cdot|||$. 
 However, it overestimates the $L^2$-norm $|||\cdot|||$ which is of order $h^2$, since the majorant
 only decreases with order $h$.
 Following the idea from \cite{Wolfmayr2016}
 a weighted
$H^1$-norm is introduced 
decreasing with the same order as the majorant.
For that, we define the norm
$ |||u-v|||_1^2 := \frac{1}{2} \|y(u) - y(v)\|^2
 + \frac{\lambda \underline{ \mu_1}^2}{2 C_F^2}  |y(u) - y(v)|_{1,\frac{1}{2}}^2$.
\begin{theorem}
 \label{theorem:DiffCombNormH1}
 The following estimate: 
 \begin{align}
  \label{equation:DiffCombNormH1}
  |||u-v|||_1^2 \leq 
  \mathcal{J}(y(v),v) - \mathcal{J}(y(u),u)
  + \frac{3 \lambda}{4 C_F^2}\left(C_F \|\mathcal{R}_1(\eta,\boldsymbol{\tau},v)\| + \|\mathcal{R}_2(\eta,\boldsymbol{\tau})\|\right)^2
 \end{align}
 is valid
 for an arbitrary control function $v \in L^2(Q)$ and
  with
   $\mathcal{R}_1(\eta,\boldsymbol{\tau},v)$ and
$\mathcal{R}_2(\eta,\boldsymbol{\tau})$ defined as in
\eqref{def:R1R2}.
\begin{proof}
Let the parameter $\delta > 0$ be 
arbitrary but fixed. 
We add and subtract
 $\eta$, apply triangle inequality on 
 $\frac{\underline{\mu_1}^2}{C_F^2 \delta}  |y(u) - y(v)|_{1,\frac{1}{2}}^2$ and
 obtain 
 \begin{align*}
 \frac{\underline{\mu_1}^2}{C_F^2 \delta}  |y(u) - y(v)|_{1,\frac{1}{2}}^2
 \leq \frac{\underline{\mu_1}^2}{2 C_F^2 \delta}
 \left( |y(u) - \eta|_{1,\frac{1}{2}}^2 +  |y(v) - \eta|_{1,\frac{1}{2}}^2\right),
 \end{align*}
where we further add and subtract $v$.
Then we apply triangle inequality two times leading to
  \begin{align*}
 \frac{\underline{\mu_1}^2}{C_F^2 \delta}  |y(u) - y(v)|_{1,\frac{1}{2}}^2
  \leq \frac{1}{2 C_F^2 \delta}
 \bigg(&\big(\|\boldsymbol{\tau}-\nu \nabla \eta\| 
 + C_F \|\sigma \partial_t \eta - \nabla \cdot \boldsymbol{\tau} - v\| + C_F \|u -v \| 
 \big)^2 \\
 &+ \left(\|\boldsymbol{\tau}-\nu \nabla \eta\|
 + C_F \|\sigma \partial_t \eta - \nabla \cdot \boldsymbol{\tau} - v\|\right)^2\bigg) \\
   \leq \frac{3}{2 C_F^2 \delta}
 (\|&\boldsymbol{\tau}-\nu \nabla \eta\| 
 + C_F \|\sigma \partial_t \eta - \nabla \cdot \boldsymbol{\tau} - v\|)^2
 + \frac{1}{\delta} \|u -v \|^2.
 \end{align*}
Together with
\eqref{equation:DiffCombNorm}
this yields
   \begin{align*}
 &|||u-v|||^2 
 + \frac{\underline{\mu_1}^2}{C_F^2 \delta}  |y(u) - y(v)|_{1,\frac{1}{2}}^2
 - \frac{1}{\delta} \|u -v \|^2 
 = \frac{1}{2} \|y(u) - y(v)\|^2
 + \frac{\underline{\mu_1}^2}{C_F^2 \delta}  |y(u) - y(v)|_{1,\frac{1}{2}}^2 \\
 &+ \left(\frac{\lambda}{2} - \frac{1}{\delta} \right) \|u-v\|^2 
 \leq \mathcal{J}(y(v),v) - \mathcal{J}(y(u),u) 
 + \frac{3}{2 C_F^2 \delta}
 \left(C_F \|\mathcal{R}_1(\eta,\boldsymbol{\tau},v)\| + \|\mathcal{R}_2(\eta,\boldsymbol{\tau})\|\right)^2.
 \end{align*} 
 We see that 
 the choice $\delta = 2/\lambda$
 finally provides estimate \eqref{equation:DiffCombNormH1}.
\end{proof}
\end{theorem}
This theorem directly leads to the following two results
presented in Propositions \ref{theorem:majorantCombNorm1}
and \ref{proposition:overallmajorantH11}.
\begin{proposition}
\label{theorem:majorantCombNorm1} 
The following 
error majorant for any control 
$v \in L^2(Q)$ is obtained:
\begin{align}
\begin{aligned}
 \label{definition:majorantCombinedNorm1}
  |||u-v|||^2_1 \leq  \mathcal{M}_1^\oplus(\alpha,\beta;\eta,\zeta,\boldsymbol{\tau},\boldsymbol{\rho},v) 
  := & \, \mathcal{J}^\oplus(\alpha,\beta;\eta,\boldsymbol{\tau},v)
  - \mathcal{J}^\ominus(\eta,{\zeta},\boldsymbol{\tau},\boldsymbol{\rho}) \\
  &+ \frac{3 \lambda}{4 C_F^2}\left(C_F \|\mathcal{R}_1(\eta,\boldsymbol{\tau},v)\| + \|\mathcal{R}_2(\eta,\boldsymbol{\tau})\|\right)^2
  \end{aligned}
\end{align}
with
\begin{align*}
 \mathcal{M}_1^\oplus(\alpha,\beta&;\eta,\zeta,\boldsymbol{\tau},\boldsymbol{\rho},v) =
   \frac{\alpha}{2} \|\eta - y_d\|^2
     + \gamma ( 
     \|\mathcal{R}_2(\eta,\boldsymbol{\tau})\|^2 
     + \frac{C_F^2}{\beta}
      \|\mathcal{R}_1(\eta,\boldsymbol{\tau},v)\|^2 )
     + \frac{\lambda}{2} \|v\|^2 - \frac{1}{2 \lambda} \|\zeta\|^2 \\
   & + \intQT \left(\nu \nabla \eta \cdot \nabla \zeta
      + \sigma \partial_t \eta \, \zeta
      + \lambda^{-1} \zeta^2 \right) \, \dxdt + \frac{C_F^2}{\underline{\mu_1}^2 \lambda}
  \left(C_F \|\mathcal{R}_3(\zeta,\boldsymbol{\rho},\eta)\| 
  + \|\mathcal{R}_4(\zeta,\boldsymbol{\rho})\| \right)^2 \\
 &+ \frac{1}{\underline{\mu_1}} (C_F 
 \|\mathcal{R}_1(\eta,\boldsymbol{\tau},-\lambda^{-1} \zeta)\| 
      + \|\mathcal{R}_2(\eta,\boldsymbol{\tau})\|) 
  \left(C_F \|\mathcal{R}_3(\zeta,\boldsymbol{\rho},\eta)\| 
  + \|\mathcal{R}_4(\zeta,\boldsymbol{\rho})\| \right) \\
  &+ \frac{3 \lambda}{4 C_F^2}\left(C_F \|\mathcal{R}_1(\eta,\boldsymbol{\tau},v)\| + \|\mathcal{R}_2(\eta,\boldsymbol{\tau})\|\right)^2
\end{align*}
where $\underline{\mu_{1}} = \frac{1}{\sqrt{2}} \min\{\underline{\nu},\underline{\sigma}\}$, 
$\alpha,\beta > 0$
as well as 
arbitrary
$\eta, \zeta \in H^{1,1}_{0,per}(Q)$ and
$\boldsymbol{\tau}, \boldsymbol{\rho}
 \in H(\emph{div},Q)$.
\end{proposition}
\begin{proposition}
\label{proposition:overallmajorantH11}
 The infimum of the majorant
 $\mathcal{M}_1^\oplus(\alpha,\beta;\eta,\zeta,\boldsymbol{\tau},\boldsymbol{\rho},v)$ 
 in \eqref{definition:majorantCombinedNorm1} 
 is attained 
 for the 
 minimum 
 of
 the optimization 
 problem I 
 being equivalent to solving
 the optimality system \eqref{problem:KKTSysSTVFAPost1}--\eqref{problem:KKTSysSTVFAPost2}
 $(v=u$, $\eta = y(u)$, ${\zeta}= p(u) = - \lambda u$,
 $\boldsymbol{\tau} = \nu \nabla y(u)$,
 $\boldsymbol{\rho} = \nu \nabla p(u))$ as follows
\begin{align*}
    \inf_{\substack{\eta,\zeta \in H^{1,1}_{0,per}(Q), 
    \boldsymbol{\tau}, \boldsymbol{\rho} \in H(\emph{div},Q), \\
     v \in L^2(Q), \alpha,\beta > 0}}
  \mathcal{M}_1^\oplus(\alpha,\beta;\eta,\zeta,\boldsymbol{\tau},\boldsymbol{\rho},v) = 0,
\end{align*}
and $\alpha$ going to zero.
\end{proposition}

\section{Two-sided bounds for optimization problem II}
\label{Sec:TwoSidedBoundsII}

In this section, we analogously derive the majorants and minorants for 
the second cost functional, however, skipping details which are similar to the derivation
in the case of optimization problem I. 

\subsection{Majorant for cost functional \eqref{equation:minfunc2}}
\label{SSec:MajorantII}

Adding and subtracting $\nabla \eta$ in the cost functional $\tilde{\mathcal{J}}(y(v),v)$,
applying the triangle inequality as well as using the estimate
  \begin{align*}
   \|\nabla y(v) - \nabla \eta\|^2
    \leq |y(v) - \eta|_{1,\frac{1}{2}}^2 = \|\nabla y(v) - \nabla \eta\|^2
    + \|\partial_t^{1/2} y(v) - \partial_t^{1/2} \eta\|^2,
  \end{align*}
we conclude that 
  \begin{align*}
   \tilde{\mathcal{J}}(y(v),v)
   &\leq \frac{1}{2} \left(\|\nabla \eta - \bsg_d\| + |y(v) - \eta|_{1,\frac{1}{2}} \right)^2
     + \frac{\lambda}{2} \|v\|^2.
  \end{align*}  
  Together with \eqref{inequality:aposteriorEstimateH11/2SeminormBVP}
this leads to the estimate 
  \begin{align*}
    \tilde{\mathcal{J}}(y(v),v)
   &\leq \frac{1}{2} \left(\|\nabla \eta - \bsg_d\| 
     + \frac{1}{\underline{\mu_1}} \|\mathcal{R}_2(\eta,\boldsymbol{\tau})\|
     + \frac{C_F}{\underline{\mu_1}} \|\mathcal{R}_1(\eta,\boldsymbol{\tau},v)\| \right)^2
     + \frac{\lambda}{2} \|v\|^2,
  \end{align*}
  where again
$\underline{\mu_{1}} = \frac{1}{\sqrt{2}} \min\{\underline{\nu},\underline{\sigma}\}$
as well as $\mathcal{R}_1(\eta,\boldsymbol{\tau},v)$ and
$\mathcal{R}_2(\eta,\boldsymbol{\tau})$ are defined as in
\eqref{def:R1R2}.
Finally, introducing parameters $\alpha, \beta > 0$ and applying Young's inequality,
  we can reformulate the estimate such that the
  right-hand side is given by a quadratic functional as follows
  \begin{align*}
    \tilde{\mathcal{J}}(y(v),v)
   \leq \tilde{\mathcal{J}}^\oplus(\alpha,\beta;\eta,\boldsymbol{\tau},v) \qquad \forall \, v \in L^2(Q),
  \end{align*}
  where 
  \begin{align}
  \label{definition:majorantCostFuncII}
  \begin{aligned}
   \tilde{\mathcal{J}}^\oplus(\alpha,\beta;\eta,\boldsymbol{\tau},v) :=
   &\, \frac{1+\alpha}{2} \|\nabla \eta - \bsg_d\|^2
     + \gamma ( 
      \|\mathcal{R}_2(\eta,\boldsymbol{\tau})\|^2 
     + \frac{C_F^2}{\beta}
      \|\mathcal{R}_1(\eta,\boldsymbol{\tau},v)\|^2 )
     + \frac{\lambda}{2} \|v\|^2.
  \end{aligned}
  \end{align}
  The infimum of the majorant \eqref{definition:majorantCostFuncII} is attained for the
 minimum of the cost functional 
 \eqref{equation:minfunc2} subject to
 \eqref{equation:forwardpde:OCP1}--\eqref{equation:forwardpde:OCP3}, which is equivalent
 to the optimal value of 
 the optimality system \eqref{problem:KKTSysSTVFAPostII1}--\eqref{problem:KKTSysSTVFAPostII2}.
Analogously to \eqref{definition:majorantCostFuncInfI}, we can show that
 \begin{align}
  \label{definition:majorantCostFuncInfII}
   \inf_{\substack{\eta \in H^{1,1}_{0,per}(Q),\boldsymbol{\tau} \in H(\text{div},Q) \\
   v \in L^2(Q), \alpha, \beta > 0}}
   \tilde{\mathcal{J}}^\oplus(\alpha,\beta;\eta,\boldsymbol{\tau},v) = \tilde{\mathcal{J}}(y(u),u),
 \end{align}
and that
  \begin{align}
  \label{estimate:majCostFuncII}
  \begin{aligned}
   \tilde{\mathcal{J}}(y(u),u)
   \leq \tilde{\mathcal{J}}^\oplus(\alpha,\beta;\eta,\boldsymbol{\tau},v) \quad
   \forall \, \eta \in H^{1,1}_{0,per}(Q), \, 
   \boldsymbol{\tau} \in H(\text{div},Q), \, 
   v \in L^2(Q), \, \alpha, \beta > 0.
  \end{aligned}
  \end{align}

\subsection{Minorant for cost functional \eqref{equation:minfunc2}}
\label{SSec:MinorantII}

Let us derive now the minorant.
For any $\eta \in H^{1,1}_{0,per}(Q)$, 
we have that
\begin{align*}
\tilde{\mathcal{J}}(y(v),v) =  \frac{1}{2} \| \nabla y - \nabla \eta \|^2
+ \intQT \left(\nabla y- \nabla \eta \right) \cdot \left(\nabla \eta - \bsg_d \right) \dxdt
+ \frac{1}{2} \| \nabla \eta - \bsg_d \|^2 + \frac{\lambda}{2}  \| v \|^2
\end{align*}
for all $v \in L^2(Q)$.
The first norm is again greater or equal to zero, together with 
the identity $v = - \lambda^{-1} p(v)$, we can estimate $\tilde{\mathcal{J}}$ from below by
\begin{align*}
\tilde{\mathcal{J}}(y(v),v) = 
 \tilde{\mathcal{J}}(y(v),- \lambda^{-1} p(v))
  \geq \frac{1}{2} \| \nabla \eta - \bsg_d \|^2
 + \frac{1}{2 \lambda} \|p\|^2 
 + \intQT \left(\nabla y- \nabla \eta \right) \cdot \left(\nabla \eta - \bsg_d \right) \dxdt.
\end{align*}
Note that Remark \ref{Rem1} applies here as well.
For $\eta \in H^{1,1}_{0,per}(Q)$,
let $p_\eta,  \tilde{p}_\eta \in H^{1,\frac{1}{2}}_{0,per}(Q)$ be the solutions to the equations
\begin{align}
 \label{equation:NOCforPeta21}
   &\int_{Q} \Big(  \nu \nabla p_\eta \cdot \nabla z
      - \sigma \partial_t^{1/2} p_\eta \, \partial_t^{1/2} z^{\perp} \Big)\dxdt
      = \int_{Q} (\nabla \eta - \bsg_d) \cdot \nabla z \, \dxdt,
      \qquad \forall z \in H^{1,\frac{1}{2}}_{0,per}(Q), \\
 \label{equation:NOCforPeta22}
    &\int_{Q} \Big( \nu \nabla \eta \cdot \nabla q
      + \sigma \partial_t^{1/2} \eta \, \partial_t^{1/2} q^{\perp}
      + \lambda^{-1} \tilde{p}_\eta \,q \Big)\,d\boldsymbol{x}\,dt = 0,  
      \qquad \forall q \in H^{1,\frac{1}{2}}_{0,per}(Q).
\end{align}
Deriving the minorant for the second minimization functional uses similar ideas now as presented
for problem I
(see Subsection \ref{SSec:MinorantI}). However, in the following we present the main steps which are important
for problem II.
Adding and subtracting $p_\eta$ together with $\frac{1}{2 \lambda} \|p - p_\eta\|^2 \geq 0$
yields
\begin{align*}
\tilde{\mathcal{J}}(&y(v),v) = 
 \tilde{\mathcal{J}}(y(v),- \lambda^{-1} p(v)) 
 \\
 &\geq \frac{1}{2} \| \nabla \eta - \bsg_d \|^2
 + \frac{1}{2 \lambda} \|p_\eta\|^2 
 + \intQT \left(\nabla y- \nabla \eta \right) \cdot \left(\nabla \eta - \bsg_d \right) \dxdt
 + \intQT \lambda^{-1} \left(p - p_\eta\right) p_\eta \, \dxdt .
\end{align*}
Appying equation \eqref{equation:NOCforPeta21}, 
identity \eqref{def:identityH01/2per} (analogously to \eqref{est:minfunc1long})
and then using equations
\eqref{equation:forwardpde:VF} and \eqref{equation:NOCforPeta22} provides 
the estimate 
\begin{align*}
 \tilde{\mathcal{J}}(y(u),u) 
\geq \, \frac{1}{2}  \| \nabla \eta - \bsg_d \|^2
 +  \frac{1}{2 \lambda} \|p_\eta\|^2
 + \intQT \lambda^{-1} (\tilde{p}_\eta - p_\eta) p_\eta \, \dxdt. 
\end{align*}
Together with introducing 
an arbitrary function $\zeta \in H^{1,1}_{0,per}(Q)$, 
following \eqref{preest},
applying Theorem \ref{theorem:estimate:adjointState}
and using \eqref{preest2}, we finally derive the estimate
\begin{align}
\label{inequality:finalminorant2}
 \tilde{\mathcal{J}}(y(u),u) \geq 
 \tilde{\mathcal{J}}^\ominus(\eta,{\zeta},\boldsymbol{\tau},\boldsymbol{\rho}) \qquad
 \forall \, \eta, \zeta \in H^{1,1}_{0,per}(Q), 
  \boldsymbol{\tau},\boldsymbol{\rho} 
  \in H(\text{div},Q)
\end{align}
with the fully computable minorant 
\begin{align}
 \label{definition:minorant2}
\begin{aligned}
  \tilde{\mathcal{J}}^\ominus(\eta,{\zeta}&,\boldsymbol{\tau},\boldsymbol{\rho})
 =  \frac{1}{2} \| \nabla \eta - \bsg_d \|^2
 +  \frac{1}{2 \lambda} \|\zeta\|^2
 - \intQT \left(\nu \nabla \eta \cdot \nabla \zeta
      + \sigma \partial_t \eta \, \zeta
      + \lambda^{-1} \zeta^2 \right) \, \dxdt \\
&- \frac{C_F^2}{\underline{\mu_{1}}^2 \lambda}
  \left(C_F \|\mathcal{R}_3(\zeta,\boldsymbol{\rho},\eta)\| 
  + \|\mathcal{R}_4(\zeta,\boldsymbol{\rho})\| \right)^2 \\
 & - \frac{1}{\underline{\mu_{1}}} (C_F \|\mathcal{R}_1(\eta,\boldsymbol{\tau},-\lambda^{-1} \zeta)\| 
      + \|\mathcal{R}_2(\eta,\boldsymbol{\tau})\|) 
  \left(C_F \|\mathcal{R}_3(\zeta,\boldsymbol{\rho},\eta)\| 
  + \|\mathcal{R}_4(\zeta,\boldsymbol{\rho})\| \right).
\end{aligned}
\end{align}
The following theorem is analogous to Theorem \ref{theorem2} and
so is its proof. Also Remark \ref{remark:convergence} can be applied for optimization problem II.
\begin{theorem} 
 The supremum 
 of the minorant $\tilde{\mathcal{J}}^\ominus$ as
 defined in \eqref{definition:minorant2} is attained for  
 the minimum 
 of the cost functional 
 \eqref{equation:minfunc2} subject to
 \eqref{equation:forwardpde:OCP1}--\eqref{equation:forwardpde:OCP3}, which is equivalent
 to the optimal value of 
 the optimality system \eqref{problem:KKTSysSTVFAPostII1}--\eqref{problem:KKTSysSTVFAPostII2}
 as follows
\begin{align}
 \label{definition:minorantSup2}
  \sup_{\substack{\eta, \zeta \in H^{1,1}_{0,per}(Q),
  \boldsymbol{\tau}, \boldsymbol{\rho} \in H(\emph{div},Q)}}
  \tilde{\mathcal{J}}^\ominus(\eta,{\zeta},\boldsymbol{\tau},\boldsymbol{\rho})
  = \tilde{\mathcal{J}}(y(u),u).
\end{align}
\end{theorem}

\subsection{Guaranteed upper bounds for the discretization errors of the control 
and the state}
\label{SSec:MajControlStateII}

We present the majorants
for the control-state 
errors 
measured in the norm
$ ||||u-v||||^2 := \frac{1}{2} \|\nabla y(u) - \nabla y(v)\|^2 + \frac{\lambda}{2} \|u-v\|^2 
 = \frac{1}{2} \|\nabla y(u) - \nabla y(v)\|^2 + \frac{1}{2\lambda} \|p(u) - p(v)\|^2$.
We obtain the identity 
$  ||||u-v||||^2 = \tilde{\mathcal{J}}(y(v),v) - \tilde{\mathcal{J}}(y(u),u)$
for an arbitrary control 
$v \in L^2(Q)$ yielding 
the 
error majorant
\begin{align}
 \label{definition:majorantCombinedNorm2}
  ||||u-v||||^2 \leq \tilde{\mathcal{M}}^\oplus(\alpha,\beta;\eta,\zeta,\boldsymbol{\tau},\boldsymbol{\rho},v)
  := \tilde{\mathcal{J}}^\oplus(\alpha,\beta;\eta,\boldsymbol{\tau},v)
  - \tilde{\mathcal{J}}^\ominus(\eta,{\zeta},\boldsymbol{\tau}&,\boldsymbol{\rho})
\end{align}
for $\alpha,\beta > 0$, arbitrary $\eta, \zeta \in H^{1,1}_{0,per}(Q)$,
$  \boldsymbol{\tau},\boldsymbol{\rho} 
  \in H(\text{div},Q)$
with
\begin{align*}
 \tilde{\mathcal{M}}^\oplus(\alpha,\beta&;\eta,\zeta,\boldsymbol{\tau},\boldsymbol{\rho},v)
 =  \, \frac{\alpha}{2} \|\nabla \eta - \bsg_d\|^2
     + \gamma ( 
     \|\mathcal{R}_2(\eta,\boldsymbol{\tau})\|^2 
     + \frac{C_F^2}{\beta}
      \|\mathcal{R}_1(\eta,\boldsymbol{\tau},v)\|^2 )
     + \frac{\lambda}{2} \|v\|^2 
     - \frac{1}{2 \lambda} \|\zeta\|^2 \\
 &+ \intQT \left(\nu \nabla \eta \cdot \nabla \zeta
      + \sigma \partial_t \eta \, \zeta
      + \lambda^{-1} \zeta^2 \right) \, \dxdt 
+ \frac{C_F^2}{\underline{\mu_{1}}^2 \lambda}
  \left(C_F \|\mathcal{R}_3(\zeta,\boldsymbol{\rho},\eta)\| 
  + \|\mathcal{R}_4(\zeta,\boldsymbol{\rho})\| \right)^2 \\
 & + \frac{1}{\underline{\mu_{1}}} (C_F \|\mathcal{R}_1(\eta,\boldsymbol{\tau},-\lambda^{-1} \zeta)\| 
      + \|\mathcal{R}_2(\eta,\boldsymbol{\tau})\|) 
  \left(C_F \|\mathcal{R}_3(\zeta,\boldsymbol{\rho},\eta)\| 
  + \|\mathcal{R}_4(\zeta,\boldsymbol{\rho})\| \right),
\end{align*}
where
$\underline{\mu_{1}} = \frac{1}{\sqrt{2}} \min\{\underline{\nu},\underline{\sigma}\}$.
 The infimum of this majorant is attained 
 for the minimum of the optimization problem II being equivalent to solving 
 the optimality system \eqref{problem:KKTSysSTVFAPostII1}--\eqref{problem:KKTSysSTVFAPostII2}
 $(v=u$, $\eta = y(u)$, ${\zeta}= p(u) = - \lambda u$,
 $\boldsymbol{\tau} = \nu \nabla y(u)$
 and $\boldsymbol{\rho} = \nu \nabla p(u))$
 as follows
\begin{align*}
    \inf_{\substack{\eta, \zeta \in H^{1,1}_{0,per}(Q), 
    \boldsymbol{\tau}, \boldsymbol{\rho}  \in H(\text{div},Q), \\
     v \in L^2(Q), \alpha,\beta > 0}}
  \tilde{\mathcal{M}}^\oplus(\alpha,\beta;\eta,\zeta,\boldsymbol{\tau},\boldsymbol{\rho} ,v) = 0.
\end{align*}
Analogously, defining
$ ||||u-v||||_1^2 := \frac{1}{2} \|\nabla y(u) - \nabla y(v)\|^2
 + \frac{\lambda \underline{\mu_1}^2}{2 C_F^2}  |y(u) - y(v)|_{1,\frac{1}{2}}^2$
we derive 
\begin{align*}
 \begin{aligned}
  ||||u-v||||^2_1 \leq \tilde{\mathcal{M}}_1^\oplus(\alpha,\beta;\eta,\zeta,\boldsymbol{\tau},\boldsymbol{\rho},v) 
    := & \, \tilde{\mathcal{J}}^\oplus(\alpha,\beta;\eta,\boldsymbol{\tau},v)
  - \tilde{\mathcal{J}}^\ominus(\eta,{\zeta},\boldsymbol{\tau},\boldsymbol{\rho}) \\
  &+ \frac{3 \lambda}{4 C_F^2}\left(C_F \|\mathcal{R}_1(\eta,\boldsymbol{\tau},v)\| + \|\mathcal{R}_2(\eta,\boldsymbol{\tau})\|\right)^2,
  \end{aligned}
\end{align*}
where now
\begin{align*}
 \tilde{\mathcal{M}}_1^\oplus(\alpha,\beta&;\eta,\zeta,\boldsymbol{\tau},\boldsymbol{\rho},v) =
   \frac{\alpha}{2} \|\nabla \eta - \bsg_d\|^2
     + \gamma ( 
     \|\mathcal{R}_2(\eta,\boldsymbol{\tau})\|^2 
     + \frac{C_F^2}{ \beta }
      \|\mathcal{R}_1(\eta,\boldsymbol{\tau},v)\|^2 )
     + \frac{\lambda}{2} \|v\|^2 - \frac{1}{2 \lambda} \|\zeta\|^2 \\
   & + \intQT \left(\nu \nabla \eta \cdot \nabla \zeta
      + \sigma \partial_t \eta \, \zeta
      + \lambda^{-1} \zeta^2 \right) \, \dxdt + \frac{C_F^2}{\underline{\mu_1}^2 \lambda}
  \left(C_F \|\mathcal{R}_3(\zeta,\boldsymbol{\rho},\eta)\| 
  + \|\mathcal{R}_4(\zeta,\boldsymbol{\rho})\| \right)^2 \\
 &+ \frac{1}{\underline{\mu_1}} (C_F 
 \|\mathcal{R}_1(\eta,\boldsymbol{\tau},-\lambda^{-1} \zeta)\| 
      + \|\mathcal{R}_2(\eta,\boldsymbol{\tau})\|) 
  \left(C_F \|\mathcal{R}_3(\zeta,\boldsymbol{\rho},\eta)\| 
  + \|\mathcal{R}_4(\zeta,\boldsymbol{\rho})\| \right) \\
  &+ \frac{3 \lambda}{4 C_F^2}\left(C_F \|\mathcal{R}_1(\eta,\boldsymbol{\tau},v)\| + \|\mathcal{R}_2(\eta,\boldsymbol{\tau})\|\right)^2.
\end{align*}
All other results similar to Propositions \ref{theorem:majorantCombNorm1}
and \ref{proposition:overallmajorantH11} follow completely.

\section{Multiharmonic finite element (MhFE) discretization}
\label{Sec:MhFEM}

The desired state $y_d$ and desired gradient $\bsg_d$ belong to $L^2(Q)$
and $[L^2(Q)]^d$, respectively. So
they can be represented as Fourier series
having Fourier coefficients from $L^2(\Omega)$.
Moreover, we assume that 
$\eta$ and $\zeta$ approximating 
the exact state $y$ and adjoint state $p$,
respectively, 
as well as the vector-valued functions
$\boldsymbol{\tau}$ and $\boldsymbol{\rho}$
are given by 
truncated Fourier series.
We also have the multiharmonic time derivative defined by
$ \partial_t \eta(\boldsymbol{x},t) = \sum_{k=1}^N \left(k \omega \, \eta_k^s(\boldsymbol{x}) \cos(k \omega t)
    - k \omega \, \eta_k^c(\boldsymbol{x}) \sin(k \omega t)\right)$
as well as the gradient and divergence by
\begin{align*}
 \nabla \eta(\boldsymbol{x},t) &= \nabla \eta_0^c(\boldsymbol{x})
 + \sum_{k=1}^N \left(\nabla \eta_k^c(\boldsymbol{x}) \, \cos(k \omega t)
 + \nabla \eta_k^s(\boldsymbol{x}) \, \sin(k \omega t)\right), \\
 \nabla \cdot \boldsymbol{\tau}(\boldsymbol{x},t) &= \nabla \cdot \boldsymbol{\tau}_0^c(\boldsymbol{x})
 + \sum_{k=1}^N \left(\nabla \cdot \boldsymbol{\tau}_k^c(\boldsymbol{x}) \, \cos(k \omega t)
 + \nabla \cdot \boldsymbol{\tau}_k^s(\boldsymbol{x}) \, \sin(k \omega t)\right).
\end{align*}
For the numerical treatment, we truncate the Fourier series expansions of all appearing functions
at an index  
$N \in \mathbb{N}$ creating Fourier series approximations of the functions.
Next, we approximate the Fourier coefficients
$ \boldsymbol{y}_k = (y_k^c, y_k^s)^T \in \mathbb{V}$ and
$\boldsymbol{p}_k = (p_k^c, p_k^s)^T \in \mathbb{V}$ 
of the unknown state and adjoint state functions
by finite element functions
$\boldsymbol{y}_{kh} = (y_{kh}^c, y_{kh}^s)^T \in \mathbb{V}_h$ and
$\boldsymbol{p}_{kh} = (p_{kh}^c, p_{kh}^s)^T \in \mathbb{V}_h$.
The finite element space $\mathbb{V}_h = V_h \times V_h \subset \mathbb{V}$ with
$V_h = \mbox{span} \{\phi_1, \dots, \phi_n\} $
and 
$\{\phi_i(\boldsymbol{x}): i=1,2,\dots,n_h \}$ 
is conforming. 
We have defined $n = n_h = \mbox{dim} V_h = O(h^{-d})$ and
 $h$ is the discretization parameter.
 The basis of the finite element space $V_h$
 on the triangulation $\mathcal{T}_h$, which is regular, consists of
piecewise linear and continuous 
elements.

\subsection{Optimization problem I}
\label{NumSSec:ModelProblemI}

The MhFE discretization 
(see also \cite{LangerWolfmayr:2013} for
more details on the multiharmonic finite element analysis for
such an optimal control problem)
yields the system of linear equations having 
a saddle point structure  
corresponding to 
the decoupled
problems \eqref{equation:MultiAnsVFBlockk1}--\eqref{equation:MultiAnsVFBlockk2} as follows
\begin{align}
 \label{equation:MultiFESysBlockk}
 \left( \begin{array}{cccc}
     M_h  &  0 & -K_{h,\nu} & k \omega M_{h,\sigma} \\
     0  &  M_h & -k \omega M_{h,\sigma} & -K_{h,\nu} \\
     -K_{h,\nu}  &  -k \omega M_{h,\sigma} & -\lambda^{-1} M_h & 0 \\
     k \omega M_{h,\sigma}  &  -K_{h,\nu} & 0 & -\lambda^{-1} M_h \end{array} \right) \left( \begin{array}{c}
     \underline{y}_k^c \\
     \underline{y}_k^s \\
     \underline{p}_k^c \\
     \underline{p}_k^s \end{array} \right) = \left( \begin{array}{c}
     M_h {\underline{y}_d^c}_k \\
     M_h {\underline{y}_d^s}_k \\
     0 \\
     0 \end{array} \right) 
     \forall k=1,\dots,N.
\end{align}
The  
functions 
$y _{kh}^c(\boldsymbol{x}) = \sum_{j=1}^n y_{k,j}^c \phi_j(\boldsymbol{x})$,
$y _{kh}^s(\boldsymbol{x}) = \sum_{j=1}^n y_{k,j}^s \phi_j(\boldsymbol{x})$,
$p _{kh}^c(\boldsymbol{x}) = \sum_{j=1}^n p_{k,j}^c \phi_j(\boldsymbol{x})$ 
and
$p _{kh}^s(\boldsymbol{x}) = \sum_{j=1}^n p_{k,j}^s \phi_j(\boldsymbol{x})$
are defined by the corresponding nodal function values 
$\underline{y}_k^c = ( y_{k,j}^c)_{j=1,\dots,n}, \,
\underline{y}_k^s = ( y_{k,j}^s)_{j=1,\dots,n}, \, 
\underline{p}_k^c = ( p_{k,j}^c)_{j=1,\dots,n}, \,
\underline{p}_k^s = ( p_{k,j}^s)_{j=1,\dots,n} \in \mathbb{R}^n$.
We have defined the stiffness matrix $K_{h,\nu}$
as well as the mass matrices $M_h$ and $M_{h,\sigma}$ by their entries
\begin{align*}
\begin{aligned}
 K_{h,\nu}^{ij} = \int_{\Omega} \nu \, \nabla \phi_i \cdot \nabla \phi_j \,d\boldsymbol{x}, \quad
 M_h^{ij} = \int_{\Omega} \phi_i \phi_j \,d\boldsymbol{x}, \quad
 M_{h,\sigma}^{ij} = \int_{\Omega} \sigma \, \phi_i \phi_j \,d\boldsymbol{x}.
\end{aligned}
\end{align*}
The right-hand side 
can be obtained by computing the vectors 
\begin{align*}
\begin{aligned}
M_h {\underline{y}_d^c}_k = \bigg( \int_{\Omega} {y_d^c}_k \phi_i \,d\boldsymbol{x}
 \bigg)_{i=1,\dots,n} \quad \mbox{and} \quad
M_h {\underline{y}_d^s}_k = \bigg( \int_{\Omega} {y_d^s}_k \phi_i \,d\boldsymbol{x}
 \bigg)_{i=1,\dots,n}.
\end{aligned}
\end{align*}
Note that for the case
$k=0$, hence for  
\eqref{equation:MultiAnsVFBlock01}--\eqref{equation:MultiAnsVFBlock02},
we obtain 
\begin{align}
 \label{equation:MultiFESysBlock0}
 \left( \begin{array}{cc}
     M_h  &  -K_{h,\nu} \\
     -K_{h,\nu}  &  - \lambda^{-1} M_h \end{array} \right) \left( \begin{array}{c}
     \underline{y}_0^c \\
     \underline{p}_0^c \end{array} \right) = \left( \begin{array}{c}
     M_h {\underline{y}_d^c}_0 \\
     0 \end{array} \right).
\end{align} 
Solving all the systems of linear equations finally lead to the contributions
for computing the MhFE approximations $y_{N h}(\boldsymbol{x},t)$ and
$p_{N h}(\boldsymbol{x},t)$ given by
$ y_{N h}(\boldsymbol{x},t) = y_{0h}^c(\boldsymbol{x}) + \sum_{k=1}^N (y_{kh}^c(\boldsymbol{x}) \cos(k \omega t)$
     $ + y_{kh}^s(\boldsymbol{x}) \sin(k \omega t))$, 
$ p_{N h}(\boldsymbol{x},t) = p_{0h}^c(\boldsymbol{x}) + \sum_{k=1}^N (p_{kh}^c(\boldsymbol{x}) \cos(k \omega t)
      + p_{kh}^s(\boldsymbol{x}) \sin(k \omega t))$.
For some proper fast solvers for these systems we refer to
\cite{KollmannKolmbauerLangerWolfmayrZulehner:2013, 
KrausWolfmayr:2013,
LangerWolfmayr:2013}. 
Both, majorant \eqref{definition:majorantCostFuncI} and minorant \eqref{definition:minorant}
of the cost functional  $\mathcal{J}$ can be computed by choosing the MhFE approximations 
$y_{N h}$, $p_{N h}$ and $u_{N h} = - \lambda^{-1} p_{N h}$
as $\eta$, $\zeta$ and $v$, respectively.
The arbitrary functions $\boldsymbol{\tau}$ and $\boldsymbol{\rho}$ can also be represented in 
form of multiharmonic functions 
$\boldsymbol{\tau}_{N h} = \boldsymbol{\tau}_{0h}^c(\boldsymbol{x}) 
+ \sum_{k=1}^N (\boldsymbol{\tau}_{kh}^c(\boldsymbol{x}) \cos(k \omega t) + \boldsymbol{\tau}_{kh}^s(\boldsymbol{x}) \sin(k \omega t))$
and 
$\boldsymbol{\rho}_{N h} = \boldsymbol{\rho}_{0h}^c(\boldsymbol{x}) 
+ \sum_{k=1}^N (\boldsymbol{\rho}_{kh}^c(\boldsymbol{x}) \cos(k \omega t) + \boldsymbol{\rho}_{kh}^s(\boldsymbol{x}) \sin(k \omega t))$.
Hence, majorant \eqref{definition:majorantCostFuncI}
and minorant \eqref{definition:minorant}
have a multiharmonic structure too. 
The linearity of the problem again yields the decoupling of the problems
introducing $\alpha_k, \beta_k > 0$ 
and resulting into
majorants $\mathcal{J}_k^\oplus$ and minorants $\mathcal{J}_k^\ominus$
 corresponding to every Fourier mode. We start with the majorants
 $\mathcal{J}_0^\oplus
 = \mathcal{J}_0^\oplus(\alpha_0,\beta_0;y_{0 h}^c,p_{0 h}^c,\boldsymbol{\tau}_{0 h}^c)$ 
 and
 $\mathcal{J}_k^\oplus
 = \mathcal{J}_k^\oplus(\alpha_k,\beta_k;\boldsymbol{y}_{k h},\boldsymbol{p}_{k h},\boldsymbol{\tau}_{k h})$ 
together with defining the parameter
$\gamma_k := ((1+\alpha_k)(1+\beta_k) C_F^2)/(2 \alpha_k \underline{\mu_1}^2)$.
We have that
    \begin{align}
    \begin{aligned}
    \label{def:Joplus0}
  \mathcal{J}_0^\oplus
  = \frac{1+\alpha_0}{2} \|y_{0 h}^c
  - {y_d}_0^c\|_{\Omega}^2 + \frac{1}{2 \lambda} \|p_{0 h}^c\|_{\Omega}^2 
  + \gamma_0 (
      \|{\mathcal{R}_2}^c_0 \|_{\Omega}^2
  + \frac{C_F^2}{\beta_0}
     \|{\mathcal{R}_1}^c_0 \|_{\Omega}^2)
\end{aligned}
\end{align}
and
   \begin{align}
    \begin{aligned}      
    \label{def:Joplusk}
  \mathcal{J}_k^\oplus
  = \frac{1+\alpha_k}{2} \|\boldsymbol{y}_{k h} 
      - {\boldsymbol{y}_d}_k \|_{\Omega}^2
	+ \frac{1}{2 \lambda}   \|\boldsymbol{p}_{k h} \|_{\Omega}^2 
      + \gamma_k (\|{\mathcal{R}_2}_k \|_{\Omega}^2
      + \frac{C_F^2}{\beta_k}
      \|{\mathcal{R}_1}_k \|_{\Omega}^2).
 \end{aligned}
\end{align}
Defining $\boldsymbol{\alpha}_{N+1} = (\alpha_0,\dots,\alpha_{N+1})^T$ and
$\boldsymbol{\beta}_N = (\beta_0,\dots,\beta_N)^T$, we can write
  the overall majorant as
 \begin{align}
  \label{def:costfuncI:multiharmonic:new}
  \begin{aligned}
  \mathcal{J}^\oplus(\boldsymbol{\alpha}_{N+1},\boldsymbol{\beta}_{N};y_{N h},p_{N h},\boldsymbol{\tau}_{N h})
  = T \,  \mathcal{J}_0^\oplus
      + \frac{T}{2} \sum_{k=1}^N 
 \mathcal{J}_k^\oplus 
+ \frac{1+\alpha_{N+1}}{2} \, \mathcal{E}_N.
 \end{aligned}
\end{align}
Here, the terms are
${\mathcal{R}_1}^c_0 
= \nabla \cdot \boldsymbol{\tau}_{0 h}^c - \lambda^{-1} p_{0 h}^c$,
${\mathcal{R}_2}^c_0 
= \boldsymbol{\tau}_{0 h}^c - \nu \nabla y_{0 h}^c$,
\begin{align*}
 {\mathcal{R}_1}_k 
  &= k \omega \, \sigma \boldsymbol{y}_{k h}^\perp
  + \text{\textbf{div}} \, \boldsymbol{\tau}_{k h} - \lambda^{-1} \boldsymbol{p}_{k h} \\
  &= ({\mathcal{R}_1}^c_k, 
      {\mathcal{R}_1}^s_k)^T 
  = (-k \omega \, \sigma y_{k h}^s 
     + \nabla \cdot \boldsymbol{\tau}_{k h}^c - \lambda^{-1} p_{k h}^c,
     k \omega \, \sigma y_{k h}^c 
     + \nabla \cdot \boldsymbol{\tau}_{k h}^s - \lambda^{-1} p_{k h}^s)^T 
\end{align*}  
and
$ {\mathcal{R}_2}_k 
  = \boldsymbol{\tau}_{k h} - \nu \nabla \boldsymbol{y}_{k h}
  = ({\mathcal{R}_2}^c_k, 
      {\mathcal{R}_2}^s_k 
      )^T
  = (\boldsymbol{\tau}_{k h}^c - \nu \nabla y_{k h}^c, 
     \boldsymbol{\tau}_{k h}^s - \nu \nabla y_{k h}^s)^T$.
The truncation's remainder term 
$\mathcal{E}_N := \|y_d - {y_d}_N\|^2 =
 \frac{T}{2} \sum_{k=N+1}^\infty \|\boldsymbol{y_d}_k\|_{\Omega}^2
 = \frac{T}{2} \sum_{k=N+1}^\infty \left(\|{y_d^c}_k\|_{\Omega}^2 + \|{y_d^s}_k\|_{\Omega}^2\right)$
can always be computed with any desired accuracy, since $y_d$ is known
(see also \cite{LangerRepinWolfmayr2015}).
The minorant \eqref{definition:minorant} can be written as
\begin{align}
 \label{definition:minorantMhFE}
\begin{aligned}
 \mathcal{J}^\ominus(y_{N h},p_{N h},\boldsymbol{\tau}_{N h},\boldsymbol{\rho}_{N h})
   = T \,  \mathcal{J}_0^\ominus 
      + \frac{T}{2} \sum_{k=1}^N 
 \mathcal{J}_k^\ominus 
+ \frac{\mathcal{E}_N}{2}, 
\end{aligned}
\end{align}
where $\mathcal{J}_0^\ominus = \mathcal{J}_0^\ominus(y_{0 h}^c,p_{0 h}^c,
  \boldsymbol{\tau}_{0 h}^c,\boldsymbol{\rho}_{0 h}^c)$ and
$\mathcal{J}_k^\ominus = \mathcal{J}_k^\ominus(\boldsymbol{y}_{k h},\boldsymbol{p}_{k h},
  \boldsymbol{\tau}_{k h},\boldsymbol{\rho}_{k h})$
are given by
   \begin{align}
    \begin{aligned}
    \label{def:Jominus0}
  \mathcal{J}_0^\ominus  &= \frac{1}{2} \|y_{0 h}^c
  - {y_d}_0^c\|_{\Omega}^2 + \frac{1}{2 \lambda} \|p_{0 h}^c\|_{\Omega}^2 
  - \intOm \left(\nu \nabla y_{0 h}^c \cdot \nabla p_{0 h}^c
      + \lambda^{-1} (p_{0 h}^c)^2 \right) \, \dx \\
      &- \frac{C_F^2}{\underline{\mu_{1}}^2 \lambda} 
  \left(C_F \|{\mathcal{R}_3}_0^c\|_{\Omega}  + \|{\mathcal{R}_4}_0^c\|_{\Omega} \right)^2 
  - \frac{1}{\underline{\mu_{1}}} (C_F 
           \|{\mathcal{R}_1}_0^c\|_{\Omega} + \|{\mathcal{R}_2}_0^c\|_{\Omega}) 
  \left(C_F \|{\mathcal{R}_3}_0^c\|_{\Omega} + \|{\mathcal{R}_4}_0^c\|_{\Omega} \right)
\end{aligned}
\end{align}
with ${\mathcal{R}_3}^c_0 
= \nabla \cdot \boldsymbol{\rho}_{0 h}^c + y_{0 h}^c - {y_d}_0^c$,
${\mathcal{R}_4}^c_0 
= \boldsymbol{\rho}_{0 h}^c - \nu \nabla p_{0 h}^c$,
and
   \begin{align}
    \begin{aligned}      
    \label{def:Jominusk}
  \mathcal{J}_k^\ominus
  &= \frac{1}{2} \|\boldsymbol{y}_{k h} 
      - {\boldsymbol{y}_d}_k \|_{\Omega}^2
	+ \frac{1}{2 \lambda}   \|\boldsymbol{p}_{k h} \|_{\Omega}^2 
	 - \intOm \left(\nu \nabla  \boldsymbol{y}_{k h} \cdot \nabla \boldsymbol{p}_{k h}
      - k \omega \, \sigma \boldsymbol{y}_{k h}^\perp \cdot \boldsymbol{p}_{k h}
      + \lambda^{-1} \boldsymbol{p}_{k h}^2 \right) \, \dx \\
      &- \frac{C_F^2}{\underline{\mu_{1}}^2 \lambda} 
  \left(C_F \|{\mathcal{R}_3}_k\|_{\Omega}  + \|{\mathcal{R}_4}_k\|_{\Omega} \right)^2 
  - \frac{1}{\underline{\mu_{1}}} (C_F 
           \|{\mathcal{R}_1}_k\|_{\Omega} + \|{\mathcal{R}_2}_k\|_{\Omega}) 
  \left(C_F \|{\mathcal{R}_3}_k\|_{\Omega} + \|{\mathcal{R}_4}_k\|_{\Omega} \right)
 \end{aligned}
\end{align}
with
$ {\mathcal{R}_3}_k 
  = k \omega \, \sigma \boldsymbol{p}_{k h}^\perp
  + \nabla \cdot \boldsymbol{\rho}_{k h} + \boldsymbol{y}_{k h} - {\boldsymbol{y}_d}_k 
  = ({\mathcal{R}_3}^c_k, 
      {\mathcal{R}_3}^s_k)^T 
  = (-k \omega \, \sigma p_{k h}^s 
     + \nabla \cdot \boldsymbol{\rho}_{k h}^c + y_{k h}^c - {y_d}_k^c,
     k \omega \, \sigma p_{k h}^c 
     + \nabla \cdot \boldsymbol{\rho}_{k h}^s + y_{k h}^s - {y_d}_k^s)^T $
and
$ {\mathcal{R}_4}_k 
  = \boldsymbol{\rho}_{k h} - \nu \nabla \boldsymbol{p}_{k h}
  = ({\mathcal{R}_4}^c_k, 
      {\mathcal{R}_4}^s_k 
      )^T
  = (\boldsymbol{\rho}_{k h}^c - \nu \nabla p_{k h}^c, 
     \boldsymbol{\rho}_{k h}^s - \nu \nabla p_{k h}^s)^T$.
 \begin{remark}
For any index $\bar{N} \in \mathbb{N}$, $\bar{N} > N$, the truncated remainder term
\begin{align*}
\mathcal{E}_{N,\bar{N}} := \frac{T}{2} \sum_{k=N+1}^{\bar{N}} \|\boldsymbol{y_d}_k\|_{\Omega}^2
\end{align*}
is a fully computable lower bound for the remainder term
$\mathcal{E}_N$. For any given $y_d \in L^2(Q)$, this provides an arbitrarily tight lower bound
for the minorant \eqref{definition:minorantMhFE}, which is in return optimal.
\end{remark}
The fluxes $\boldsymbol{\tau}_{0 h}^c$, $\boldsymbol{\rho}_{0 h}^c$ and
$\boldsymbol{\tau}_{k h}$,
$\boldsymbol{\rho}_{k h}$ for all $k=1,\dots,N$, denoted by
$ \boldsymbol{\tau}_{k h}
 = R_h^{\text{\tiny{flux}}}(\nu \nabla \boldsymbol{y}_{k h})$ and
$  \boldsymbol{\rho}_{k h}
 = R_h^{\text{\tiny{flux}}}(\nu \nabla \boldsymbol{p}_{k h})$
are reconstructed by lowest-order
Raviart-Thomas elements
mapping $L^2$-functions to 
$H(\text{div},\Omega)$,
see \cite{RaviartThomas:1977}
as well as \cite{LangerRepinWolfmayr2015, LangerRepinWolfmayr2016},
leading to
$ \boldsymbol{\tau}_{N h} 
 = R_h^{\text{\tiny{flux}}}(\nu \nabla y_{N h})$ and
$   \boldsymbol{\rho}_{N h} 
 = R_h^{\text{\tiny{flux}}}(\nu \nabla p_{N h})$.
We 
minimize 
$\mathcal{J}^\oplus$
with respect to the positive parameters
$\alpha_k$ and $\beta_k$
 leading to the optimized
$\boldsymbol{\alpha}_{N+1}$ and $\boldsymbol{\beta}_{N}$.
Finally, the multiharmonic majorant \eqref{def:costfuncI:multiharmonic:new} 
and 
minorant \eqref{definition:minorantMhFE} 
lead to 
upper and 
lower bounds for $\mathcal{J}$
which are guaranteed and 
computable.
The computation of the infimum of $\mathcal{J}^{\oplus}$  and the supremum of $\mathcal{J}^{\ominus}$ 
provide the minimum 
of $\mathcal{J}$, 
see also \cite{LangerRepinWolfmayr2016}.

\subsection{Optimization problem II}
\label{NumSSec:ModelProblemII}

For the second problem, we only summarize the main results and changes.
The MhFE discretization leads to the following discrete 
problem:
\begin{align}
 \label{equation:MultiFESysBlockk:2}
 \left( \begin{array}{cccc}
     K_h  &  0 & -K_{h,\nu} & k \omega M_{h,\sigma} \\
     0  &  K_h & -k \omega M_{h,\sigma} & -K_{h,\nu} \\
     -K_{h,\nu}  &  -k \omega M_{h,\sigma} & -\lambda^{-1} M_h & 0 \\
     k \omega M_{h,\sigma}  &  -K_{h,\nu} & 0 & -\lambda^{-1} M_h \end{array} \right) \left( \begin{array}{c}
     \underline{y}_k^c \\
     \underline{y}_k^s \\
     \underline{p}_k^c \\
     \underline{p}_k^s \end{array} \right) = \left( \begin{array}{c}
     K_h {\underline{g}_d^c}_k \\
     K_h {\underline{g}_d^s}_k \\
     0 \\
     0 \end{array} \right)
\end{align}
and 
\begin{align}
 \label{equation:MultiFESysBlock0:2}
 \left( \begin{array}{cc}
     K_h  &  -K_{h,\nu} \\
     -K_{h,\nu}  &  - \lambda^{-1} M_h \end{array} \right) \left( \begin{array}{c}
     \underline{y}_0^c \\
     \underline{p}_0^c \end{array} \right) = \left( \begin{array}{c}
     K_h {\underline{g}_d^c}_0 \\
     0 \end{array} \right).
\end{align}
Now, the right-hand side vectors are computed by 
\begin{align*}
\begin{aligned}
K_h {\underline{g}_d^c}_k = \bigg( \int_{\Omega} {\boldsymbol{g}_d}_k^c \cdot \nabla \phi_i \,d\boldsymbol{x}
 \bigg)_{i=1,\dots,n} \quad \mbox{and} \quad
K_h {\underline{g}_d^s}_k = \bigg( \int_{\Omega} {\boldsymbol{g}_d}_k^s \cdot \nabla \phi_i \,d\boldsymbol{x}
 \bigg)_{i=1,\dots,n}.
\end{aligned}
\end{align*}
We summarize the discrete
majorant \eqref{definition:majorantCostFuncII} and minorant \eqref{definition:minorant2}
of the cost functional  $\tilde{\mathcal{J}}$ 
computed by choosing the MhFE approximations for all the (arbitrary) functions.
Defining now
$\tilde{\mathcal{E}}_N 
:= \|\boldsymbol{g}_d - {\boldsymbol{g}_d}_N\|^2 =
 \frac{T}{2} \sum_{k=N+1}^\infty \|{\boldsymbol{g}_d}_k\|_{\Omega}^2
 = \frac{T}{2} \sum_{k=N+1}^\infty \left(\|{\bsg_d^c}_k\|_{\Omega}^2 + \|{\bsg_d^s}_k\|_{\Omega}^2\right)$
as truncation's remainder term, we write the
majorant \eqref{definition:majorantCostFuncII} 
as
   \begin{align}
  \label{def:costfuncII:multiharmonic:new}
  \begin{aligned}
  \tilde{\mathcal{J}}^\oplus(\boldsymbol{\alpha}_{N+1},\boldsymbol{\beta}_{N};y_{N h},p_{N h},
  \boldsymbol{\tau}_{N h})
  = T \,  \tilde{\mathcal{J}}_0^\oplus 
      + \frac{T}{2} \sum_{k=1}^N 
 \tilde{\mathcal{J}}_k^\oplus 
+ \frac{1+\alpha_{N+1}}{2} \, \tilde{\mathcal{E}}_N,
 \end{aligned}
\end{align}
where $\tilde{\mathcal{J}}_0^\oplus = \tilde{\mathcal{J}}_0^\oplus(\alpha_0,\beta_0;y_{0 h}^c,p_{0 h}^c,
  \boldsymbol{\tau}_{0 h}^c)$ and 
$\tilde{\mathcal{J}}_k^\oplus = \tilde{\mathcal{J}}_k^\oplus(\alpha_k,\beta_ks; \boldsymbol{y}_{k h},\boldsymbol{p}_{k h},
  \boldsymbol{\tau}_{k h})$
are given by
   \begin{align}
    \begin{aligned}
    \label{def:Joplus0:2}
  \tilde{\mathcal{J}}_0^\oplus 
  = \frac{1+\alpha_0}{2} \|\nabla y_{0 h}^c
  - {\bsg_d}_0^c\|_{\Omega}^2 + \frac{1}{2 \lambda} \|p_{0 h}^c\|_{\Omega}^2 
  + \gamma_0 ( 
      \|{\mathcal{R}_2}^c_0 \|_{\Omega}^2
  + \frac{C_F^2}{\beta_0}
     \|{\mathcal{R}_1}^c_0 \|_{\Omega}^2)
\end{aligned}
\end{align}
and
   \begin{align}
    \begin{aligned}      
    \label{def:Joplusk:2}
  \tilde{\mathcal{J}}_k^\oplus 
  = \frac{1+\alpha_k}{2} \|\nabla \boldsymbol{y}_{k h} 
      - {\boldsymbol{g}_d}_k \|_{\Omega}^2
	+ \frac{1}{2 \lambda}   \|\boldsymbol{p}_{k h} \|_{\Omega}^2 
      + \gamma_k ( 
      \|{\mathcal{R}_2}_k \|_{\Omega}^2
      + \frac{C_F^2}{\beta_k}
      \|{\mathcal{R}_1}_k \|_{\Omega}^2). 
 \end{aligned}
\end{align}
The minorant \eqref{definition:minorant2} can be written as
\begin{align}
 \label{definition:minorantMhFE:2}
\begin{aligned}
 \tilde{\mathcal{J}}^\ominus(y_{N h},p_{N h},\boldsymbol{\tau}_{N h},\boldsymbol{\rho}_{N h})
   = &\, T \,   \tilde{\mathcal{J}}_0^\ominus
      + \frac{T}{2} \sum_{k=1}^N 
  \tilde{\mathcal{J}}_k^\ominus
+ \frac{\tilde{\mathcal{E}}_N}{2}, 
\end{aligned}
\end{align}
where $\tilde{\mathcal{J}}_0^\ominus = \tilde{\mathcal{J}}_0^\ominus(y_{0 h}^c,p_{0 h}^c,
  \boldsymbol{\tau}_{0 h}^c,\boldsymbol{\rho}_{0 h}^c)$ and
  $\tilde{\mathcal{J}}_k^\ominus = \tilde{\mathcal{J}}_k^\ominus(\boldsymbol{y}_{k h},\boldsymbol{p}_{k h},
  \boldsymbol{\tau}_{k h},\boldsymbol{\rho}_{k h})$
are given by
   \begin{align}
    \begin{aligned}
    \label{def:Jominus0:2}
  \tilde{\mathcal{J}}_0^\ominus 
  &= \frac{1}{2} \|\nabla y_{0 h}^c
  - {\bsg_d}_0^c\|_{\Omega}^2 + \frac{1}{2 \lambda} \|p_{0 h}^c\|_{\Omega}^2 
  - \intOm \left(\nu \nabla y_{0 h}^c \cdot \nabla p_{0 h}^c
      + \lambda^{-1} (p_{0 h}^c)^2 \right) \, \dx \\
     & - \frac{C_F^2}{\underline{\mu_{1}}^2 \lambda} 
  \left(C_F \|{\mathcal{R}_3}_0^c\|_{\Omega}  + \|{\mathcal{R}_4}_0^c\|_{\Omega} \right)^2 
  - \frac{1}{\underline{\mu_{1}}} (C_F 
           \|{\mathcal{R}_1}_0^c\|_{\Omega} + \|{\mathcal{R}_2}_0^c\|_{\Omega}) 
  \left(C_F \|{\mathcal{R}_3}_0^c\|_{\Omega} + \|{\mathcal{R}_4}_0^c\|_{\Omega} \right)
\end{aligned}
\end{align}
and
   \begin{align}
    \begin{aligned}      
    \label{def:Jominusk:2}
& \tilde{\mathcal{J}}_k^\ominus 
  = \frac{1}{2} \|\nabla \boldsymbol{y}_{k h} 
      - {\boldsymbol{g}_d}_k \|_{\Omega}^2
	+ \frac{1}{2 \lambda}   \|\boldsymbol{p}_{k h} \|_{\Omega}^2 
	- \intOm \left(\nu \nabla  \boldsymbol{y}_{k h} \cdot \nabla \boldsymbol{p}_{k h}
      - k \omega \, \sigma \boldsymbol{y}_{k h}^\perp \cdot \boldsymbol{p}_{k h}
      + \lambda^{-1} \boldsymbol{p}_{k h}^2 \right) \, \dx \\
      &- \frac{C_F^2}{\underline{\mu_{1}}^2 \lambda} 
  \left(C_F \|{\mathcal{R}_3}_k\|_{\Omega}  + \|{\mathcal{R}_4}_k\|_{\Omega} \right)^2 
  - \frac{1}{\underline{\mu_{1}}} (C_F 
           \|{\mathcal{R}_1}_k\|_{\Omega} + \|{\mathcal{R}_2}_k\|_{\Omega}) 
  \left(C_F \|{\mathcal{R}_3}_k\|_{\Omega} + \|{\mathcal{R}_4}_k\|_{\Omega} \right).
 \end{aligned}
\end{align}

\section{Robust preconditioners for the minimal residual method}
\label{Sec:Preconditioning}

The saddle point systems \eqref{equation:MultiFESysBlockk},
\eqref{equation:MultiFESysBlock0}, 
\eqref{equation:MultiFESysBlockk:2} and \eqref{equation:MultiFESysBlock0:2}
can be solved by using 
the preconditioned MinRes method, see  
\cite{PaigeSaunders:1975}.
A convergence result for the preconditioned MinRes method can be found in 
\cite{Greenbaum:1997} stating that the convergence rate of the preconditioned
MinRes method depends on the condition number of the preconditioned system.
The derivation of preconditioners for problems  
\eqref{equation:MultiFESysBlockk} for $k=1,\dots,N$
and \eqref{equation:MultiFESysBlock0} for $k=0$ have already been 
presented and discussed in
\cite{KollmannKolmbauerLangerWolfmayrZulehner:2013, LangerWolfmayr:2013}
given by
\begin{align}
\label{equation:preconditionerSigmaPWConst:OCP}
 \mathcal{P}_k = \text{diag}(D_k,D_k,\lambda^{-1} D_k,\lambda^{-1} D_k) \quad \text{and}
\quad  \mathcal{P}_0 = \text{diag}(D_0,\lambda^{-1} D_0),
\end{align}
respectively,
where 
$D_k = \sqrt{\lambda} K_{h,\nu} + k \omega \sqrt{\lambda} M_{h,\sigma} + M_h$
and
$D_0 = M_h + \sqrt{\lambda} K_{h,\nu}$.
In \cite{KollmannKolmbauerLangerWolfmayrZulehner:2013},
preconditioners are derived following 
the technique in \cite{Zulehner2011} based on operator interpolation
theory. 

In this section, we present new
robust preconditioners for the problem matrices in
\eqref{equation:MultiFESysBlockk:2} for $k=1,\dots,N$
and in \eqref{equation:MultiFESysBlock0:2} for $k=0$
in order to solve
optimization problem II.
Here, we assume that $\sigma$ and $\nu$ are constant, which we also choose 
in the numerical results
in Section \ref{Sec:NumericalResults}. Hence,
$M_{h,\sigma} = \sigma M_h$ and $K_{h,\nu} = \nu K_h$.
The
block-diagonal preconditioners are practically 
implemented by
the version of the 
algebraic multilevel iteration (AMLI) method from
\cite{Kraus:2012}. 
The AMLI preconditioned MinRes solver
is robust and of optimal complexity which is proved in \cite{KrausWolfmayr:2013}.
This can be also observed in the numerical results' section.
Let us consider the saddle point system \eqref{equation:MultiFESysBlockk:2} for $k \geq N$.
The derivation of preconditioners 
for system \eqref{equation:MultiFESysBlock0:2} 
in case of $k = 0$ is completely analogous.
We refer the reader to
\cite{KollmannKolmbauerLangerWolfmayrZulehner:2013} on
details for the derivation of the preconditioners based
on Schur complements.

Defining the 
matrices and vectors
\begin{align}
\label{defABC}
 A = \text{diag}(K_h,K_h),
\qquad
 B = \left( \begin{array}{cc}
     -\nu K_h  &  -k \omega \sigma M_h \\
     k \omega \sigma M_h & - \nu K_h \end{array} \right), 
\qquad
     C = \text{diag}(\lambda^{-1}  M_h, \lambda^{-1}  M_h),
\end{align}
$ \underline{f} = ( K_h {\underline{g}_d^c}_k,
     K_h {\underline{g}_d^s}_k)^T$,
$ \underline{y} = (\underline{y}_k^c, \underline{y}_k^s)^T$
 and
 $\underline{p} = (\underline{p}_k^c, \underline{p}_k^s)^T$
leads to the following problem structure
\begin{align}
 \label{equation:SaddlePointProblem}
  \left( \begin{array}{cc}
     A  &  B^T \\
     B  &  -C \end{array} \right) 
     \left( \begin{array}{c}
     \underline{y} \\
     \underline{p} \end{array} \right) = 
     \left( \begin{array}{c}
     \underline{f} \\
     0 \end{array} \right) 
\end{align}
with the symmetric and positive definite matrices $A$ and $C$.
We define the negative Schur complements 
$S = C + B A^{-1} B^T$ and $
  R = A + B^T C^{-1} B$
yielding the preconditioners
for $k \geq N$ as follows
 \begin{align}
 \label{definition:precondP0P1}
  \tilde{\mathcal{P}} = \text{diag}(A,S) 
\qquad
\mbox{ and } \qquad
\tilde{\mathcal{Q}} =  \text{diag}(R,C).
 \end{align}
The negative Schur complements are given by
\begin{align}
\label{defS}
 S =  \text{diag}(\nu K_h + \lambda^{-1} M_h + k^2 \omega^2 \sigma^2  M_h K_h^{-1} M_h,
 \nu K_h + \lambda^{-1} M_h + k^2 \omega^2 \sigma^2  M_h K_h^{-1} M_h) 
\end{align}
and
\begin{align}
\label{defR}
 R = \text{diag}(K_h + k^2 \omega^2 \sigma^2 \lambda M_h + \nu^2 \lambda K_h M_h^{-1} K_h,
 K_h + k^2 \omega^2 \sigma^2 \lambda M_h + \nu^2 \lambda K_h M_h^{-1} K_h).
\end{align}
Let us define
\begin{align}
\label{defDkS}
\tilde{D}_k^S = \nu K_h + \lambda^{-1} M_h + k^2 \omega^2 \sigma^2  M_h K_h^{-1} M_h
\end{align}
and
\begin{align}
\label{defDkR}
\tilde{D}_k^R = K_h + k^2 \omega^2 \sigma^2 \lambda M_h + \nu^2 \lambda K_h M_h^{-1} K_h.
\end{align}
Then $S$ and $R$ can be written 
as $S = \text{diag}(\tilde{D}_k^S,\tilde{D}_k^S)$ 
and $R = \text{diag}(\tilde{D}_k^R,\tilde{D}_k^R)$, respectively.
\begin{remark}
Both Schur complement preconditioners $\tilde{\mathcal{P}}$ and
$\tilde{\mathcal{Q}}$ in \eqref{definition:precondP0P1}
can be chosen for computations of optimization problem II
leading to fast and robust convergence rates, see  
 \cite{Kuznetsov:1995} and 
 \cite{MurphyGolubWathen:2000}.
\end{remark}
Inserting now $A$ and $C$ from \eqref{defABC}
and $S$ and $R$ from \eqref{defS} and \eqref{defR}, respectively,
the 
Schur complement preconditioners in the form of $\tilde{\mathcal{P}}$
as presented in \eqref{definition:precondP0P1}
in the case of $k \geq N$ as well as for $k=0$ are 
given by
\begin{align}
\label{equation:preconditionerSigmaPWConst:OCP2}
 \tilde{\mathcal{P}}_k = \text{diag}(K_{h},K_{h},\tilde{D}_k^S,\tilde{D}_k^S) \quad \text{and} \quad
 \tilde{\mathcal{P}}_0 = \text{diag}(K_{h},\nu K_{h}+\lambda^{-1} M_{h})
\end{align}
Analogously in the form of $\tilde{\mathcal{Q}}$, they are given by
\begin{align}
\label{equation:preconditionerSigmaPWConst:OCP2b}
 \tilde{\mathcal{Q}}_k = \text{diag}(\tilde{D}_k^R,\tilde{D}_k^R,\lambda^{-1}  M_h, \lambda^{-1}  M_h) \quad \text{and} \quad
 \tilde{\mathcal{Q}}_0 = \text{diag}(\lambda^{-1} M_{h}, K_h + \nu^2 \lambda K_h M_h^{-1} K_h),
\end{align}
for $k \geq N$ and $k=0$
for the saddle point systems \eqref{equation:MultiFESysBlockk:2} and \eqref{equation:MultiFESysBlock0:2}, respectively.
We refer again to \cite{KollmannKolmbauerLangerWolfmayrZulehner:2013}
for further details on the preconditioners' derivation.
For the numerical experiments of this work, we simplify the preconditioner such that
$\tilde{ \tilde{\mathcal{P}}}_k = \text{diag}(\tilde{\tilde{D}}_k^S,\tilde{\tilde{D}}_k^S,\lambda^{-1} \tilde{\tilde{D}}_k^S,\lambda^{-1} \tilde{\tilde{D}}_k^S)$
and $\tilde{ \tilde{\mathcal{P}}}_0 = \text{diag}(\tilde{\tilde{D}}_0^S,\lambda^{-1} \tilde{\tilde{D}}_0^S)$, where 
$\tilde{\tilde{D}}_k^S= \nu \sqrt{\lambda} K_{h} + (1 + k \omega \sigma \sqrt{\lambda}) M_{h}$
and
$\tilde{\tilde{D}}_0^S = M_h + \nu \sqrt{\lambda} K_{h}$
in a similar form as \eqref{equation:preconditionerSigmaPWConst:OCP}
in order to apply the AMLI preconditioned MinRes method analogously as for optimization problem I.

\section{Numerical results}
\label{Sec:NumericalResults}

In this section, we present numerical results performed in C\texttt{++} for computing
the minorants and majorants of the two optimal control problems with 
cost functionals \eqref{equation:minfunc1} and \eqref{equation:minfunc2} and for different cases of given data
denoted by Examples 1--6.
For the numerical results on the majorants of optimization problem I, we refer to
\cite{LangerRepinWolfmayr2016} and for the numerical analysis including convergence results,
to
\cite{KollmannKolmbauerLangerWolfmayrZulehner:2013, LangerWolfmayr:2013}.
However, all the numerical experiments on the majorants
of optimization problem II as well as all on minorants for both problems I and II are new.
We perform numerical experiments for the same three cases (applied on the problem data) as for problem I but 
in Examples 4--6 they are applied on the desired gradients.

The 
domain 
$\Omega$ is the two-dimensional unit square $\Omega = (0,1)^2$
with Friedrichs' constant $C_F = 1/(\sqrt{2}\pi)$.
The triangulation of the domain is performed regular leading to a uniform grid.
The finite elements are chosen as described in Section \ref{Sec:MhFEM}.
The coefficients are chosen to be $\nu = \sigma = 1$.
In Examples 1, 2, 4 and 5,  
the cost parameter is chosen to be $\lambda = 0.1$,
$T=2 \pi/\omega$ and $\omega = 1$.
In Examples 3 and 6, 
$\lambda = 0.01$, $T=1$ and $\omega = 2 \pi$.
The MhFE approximations
for $\eta$,
$\zeta$ and $\boldsymbol{\tau}, \boldsymbol{\rho}$
are chosen as well as the fluxes are reconstructed by 
$RT^0$-extensions (lowest-order standard Raviart-Thomas) 
resulting in averaged fluxes which are now from $H(\text{div},\Omega)$,
see also \cite{LangerRepinWolfmayr2015,LangerRepinWolfmayr2016} for further details.
The grid sizes range between
16$\times$16 and 256$\times$256 as well as
512$\times$512 to obtain finer grid solutions for Examples 3 and 6 as a reference for the exact solution.
The preconditioned MinRes iteration was stopped after 8 iteration steps in all
computations using the AMLI preconditioner 
with 4 inner iterations.
The numerical experiments for Examples 1, 2, 4 and 5
were performed on a laptop with
Intel(R) Core(TM) i5-6267U CPU @ 2.90GHz processor and 16 GB 2133 MHz LPDDR3 memory.
The numerical experiments for Examples 3 and 6 were performed on
a CPU server with a Tumbleweed distribution having 64 cores and 1 terabyte memory
in order to provide enough memory for computing
the finer grid solutions in addition.
All computational times
$t^{\text{sec}}$ are presented in seconds and
include the CPU times also  
needed to derive the majorants and minorants. However, we want to highlight that these times are
much smaller compared to the rest.
The computational times of Examples 3 and 6 exclude the computation
of the solution on the finer grid (512$\times$512).

\subsection{Numerical results for optimization problem I}
\label{SSec:NumericalResultsI}

\subsubsection{Example 1}
\label{SSSec:Ex1}

The desired state is a time-periodic and time-analytic function
\begin{align}
\label{ydEx1}
 y_d(\boldsymbol{x},t) = e^t \sin(t) \, 0.1 \left(\left(12+4 \pi^4\right) \sin^2(t)-6 \cos(t) (\cos(t)-\sin(t))\right)
			 \sin (x_1 \pi) \sin (x_2 \pi),
\end{align}
which is however not time-harmonic. 
Note that the exact state function for this example is 
\begin{align}
\label{eqn:exactSol:ex1}
y(\boldsymbol{x},t) = e^t \sin(t)^3 \sin (x_1 \pi) \sin (x_2 \pi). 
\end{align}
The truncation index for the multiharmonic approximations is chosen as
$N = 8$ here.
Table \ref{tab:Ex1:0} presents for different grid sizes  
CPU times $t^{\text{sec}}$, 
values for the majorants $\mathcal{J}^{\oplus}_0$ and minorants $\mathcal{J}^{\ominus}_0$
as defined in \eqref{def:Joplus0} and \eqref{def:Jominus0}, 
and
efficiency indices
$I_{\text{eff}}^{\mathcal{J}^\oplus_0} = \mathcal{J}^{\oplus}_0 / \mathcal{J}_0$,
$I_{\text{eff}}^{\mathcal{J}^\ominus_0} = \mathcal{J}^{\ominus}_0 / \mathcal{J}_0$ and
$I_{\text{eff}}^{\mathcal{J},0} = \mathcal{J}^{\oplus}_0 / \mathcal{J}^{\ominus}_0$.
Here, $\mathcal{J}_0 = \mathcal{J}_0(y_0^c,u_0^c)
  = \frac{1}{2} \|y_0^c - {y_d^c}_0\|_{\Omega}^2
     + \frac{\lambda}{2} \|u_0^c\|_{\Omega}^2$
 as introduced in Subsection \ref{SSec:ModelProblemI}.
In Table~\ref{tab:Ex1:1}, 
the numerical results for the Fourier mode 
$k=1$ are presented including
$\mathcal{J}^{\oplus}_k$, $\mathcal{J}^{\ominus}_k$
as defined in \eqref{def:Joplusk} and \eqref{def:Jominusk} and
the  
corresponding efficiency indices
$I_{\text{eff}}^{\mathcal{J}^\oplus_k} = \mathcal{J}^{\oplus}_k / \mathcal{J}_k$,
$I_{\text{eff}}^{\mathcal{J}^\ominus_k} = \mathcal{J}^{\ominus}_k / \mathcal{J}_k$ and
$I_{\text{eff}}^{\mathcal{J},k} = \mathcal{J}^{\oplus}_k / \mathcal{J}^{\ominus}_k$.
Moreover, we present the efficiency indices for $\mathcal{M}_1^{\oplus}$ 
given for the modes by
\begin{align*}
   I_{\text{eff}}^{\mathcal{M}_1,0} =  
  \sqrt{\frac{\mathcal{M}_{1,0}^\oplus(\alpha_0,\beta_0;y_{0 h}^c,p_{0 h}^c,
  \boldsymbol{\tau}_{0 h}^c,\boldsymbol{\rho}_{0 h}^c)}{|||y_{0}^c - y_{0 h}^c|||_{1,0}^2}}
\quad \text{ and } \quad
  I_{\text{eff}}^{\mathcal{M}_1,k} =  
  \sqrt{\frac{\mathcal{M}_{1,k}^\oplus(\alpha_k,\beta_k;\boldsymbol{y}_{k h},\boldsymbol{p}_{k h},
  \boldsymbol{\tau}_{k h},\boldsymbol{\rho}_{k h})}{|||\boldsymbol{y}_{k} - \boldsymbol{y}_{k h}|||_{1,k}^2}}.
\end{align*}
The error norms for the modes are given by
\begin{align*}
|||y_{0}^c - y_{0 h}^c|||_{1,0}^2 =
\frac{1}{2} \|y_{0}^c - y_{0 h}^c\|_{\Omega}^2
 + \frac{\lambda \underline{\mu_1}^2}{2 C_F^2}  \|\nabla y_0^c - \nabla y_{0 h}^c\|_{\Omega}^2
 \end{align*}
and
\begin{align*}
|||\boldsymbol{y}_{k} - \boldsymbol{y}_{k h}|||_{1,k}^2 =
\left(\frac{1}{2}  + \frac{k\omega \lambda \underline{\mu_1}^2}{2 C_F^2}   \right)
  \|\boldsymbol{y}_{k} - \boldsymbol{y}_{k h}\|_{\Omega}^2
 + \frac{\lambda \underline{\mu_1}^2}{2 C_F^2}   \|\nabla \boldsymbol{y}_{k} - \nabla \boldsymbol{y}_{k h}\|_{\Omega}^2
 \end{align*}
leading the representation
   \begin{align}
  \label{definition:combinedNormH1withk}
  \begin{aligned}
 |||u-v|||_1^2 
  = &\, T \,  |||y_{0}^c - y_{0 h}^c|||_{1,0}^2 
      + \frac{T}{2} \sum_{k=1}^N  |||\boldsymbol{y}_{k} - \boldsymbol{y}_{k h}|||_{1,k}^2 
+ \mathcal{F}_N
 \end{aligned}
\end{align}  
with the remainder term
$ \mathcal{F}_N :=  \frac{T}{2} \sum_{k={N+1}}^\infty  |||\boldsymbol{y}_{k}|||_{1,k}^2$. 
For the numerical experiments, we can estimate the efficiency index
for $\mathcal{M}_1^{\oplus}$ 
from above by estimating \eqref{definition:combinedNormH1withk}
from below ignoring the remainder term $\mathcal{F}_N$ 
leading to the overall efficiency index for $\mathcal{M}_1^{\oplus}$ 
\begin{align}
\label{IeffM1New}
  I_{\text{eff}}^{\mathcal{M}_1} =  
  \sqrt{\frac{\mathcal{M}_1^\oplus(\alpha,\beta;\eta,\zeta,\boldsymbol{\tau},\boldsymbol{\rho},v)}{T \,  |||y_{0}^c - y_{0 h}^c|||_{1,0}^2  + \frac{T}{2} \sum_{k=1}^N  |||\boldsymbol{y}_{k} - \boldsymbol{y}_{k h}|||_{1,k}^2
}}.
\end{align}
The corresponding majorants 
$\mathcal{M}_{1,0}^\oplus = \mathcal{M}_{1,0}^\oplus(\alpha_0,\beta_0;y_{0 h}^c,p_{0 h}^c,
  \boldsymbol{\tau}_{0 h}^c,\boldsymbol{\rho}_{0 h}^c) $
and $\mathcal{M}_{1,k}^\oplus =
  \mathcal{M}_{1,k}^\oplus(\alpha_k,\beta_k;$ $\boldsymbol{y}_{k h},\boldsymbol{p}_{k h},
  \boldsymbol{\tau}_{k h},\boldsymbol{\rho}_{k h})$
are given by
   \begin{align*}
    \mathcal{M}_{1,0}^\oplus
  &= \mathcal{J}_0^\oplus 
  - \mathcal{J}_0^\ominus 
  + \frac{3 \lambda}{4 C_F^2}\left(C_F \|{\mathcal{R}_1}^c_0 \|_{\Omega} +  \|{\mathcal{R}_2}^c_0 \|_{\Omega} \right)^2 
  = \frac{\alpha_0}{2} \|y_{0 h}^c  - {y_d}_0^c\|_{\Omega}^2 
  + \gamma_0 
      \|{\mathcal{R}_2}^c_0 \|_{\Omega}^2 \\
  &+ \frac{\gamma_0 C_F^2}{\beta_0}
     \|{\mathcal{R}_1}^c_0 \|_{\Omega}^2  
     + \intOm \left(\nu \nabla y_{0 h}^c \cdot \nabla p_{0 h}^c
      + \lambda^{-1} (p_{0 h}^c)^2 \right) \, \dx 
      + \frac{C_F^2}{\underline{\mu_{1}}^2 \lambda} 
  \left(C_F \|{\mathcal{R}_3}_0^c\|_{\Omega}  + \|{\mathcal{R}_4}_0^c\|_{\Omega} \right)^2 \\
  &+ \frac{1}{\underline{\mu_{1}}} (C_F 
           \|{\mathcal{R}_1}_0^c\|_{\Omega} + \|{\mathcal{R}_2}_0^c\|_{\Omega}) 
  \left(C_F \|{\mathcal{R}_3}_0^c\|_{\Omega} + \|{\mathcal{R}_4}_0^c\|_{\Omega} \right)  
    + \frac{3 \lambda}{4 C_F^2}\left(C_F \|{\mathcal{R}_1}^c_0 \|_{\Omega} +  \|{\mathcal{R}_2}^c_0 \|_{\Omega} \right)^2
\end{align*}
and also 
$\mathcal{M}_{1,k}^\oplus 
 =  \mathcal{J}_k^\oplus 
  -  \mathcal{J}_k^\ominus 
  + \frac{3 \lambda}{4 C_F^2}\left(C_F 
           \|{\mathcal{R}_1}_k\|_{\Omega} + \|{\mathcal{R}_2}_k\|_{\Omega} \right)^2$.
Table~\ref{tab:Ex1:global}
sums up the numerical results for Example 1 by presenting
the minorants, majorants and efficiency indices 
on a grid of size $256\times256$
for all $k$ up to $N = 4$
and then for $k = 6$ and $k = 8$
(since their results were similar as for $k = 5$ and $k = 7$).
\begin{table}[!ht]
\begin{center}
\begin{tabular}{|c|c|cc|cc|c|c|}
  \hline
   grid & $t^{\text{sec}}$ 
   & $\mathcal{J}^{\ominus}_0$ 
   & $I_{\text{eff}}^{\mathcal{J}^\ominus_0}$
   & $\mathcal{J}^{\oplus}_0$ 
   & $I_{\text{eff}}^{\mathcal{J}^\oplus_0}$ 
   & $I_{\text{eff}}^{\mathcal{J},0}$ 
   & $I_{\text{eff}}^{\mathcal{M}_1,0}$ \\
  \hline
   $16   \times   16$   &  0.02  & 1.13e+05 & 0.90 & 1.26e+05 & 1.01 & 1.12 & 1.53  \\ 
   $32   \times   32$   &  0.07  & 1.14e+05 & 0.90 & 1.27e+05 & 1.00 & 1.11 & 1.47  \\ 
   $64   \times   64$   &  0.24  & 1.14e+05 & 0.90 & 1.27e+05 & 1.00 & 1.11 & 1.44  \\   
   $128 \times 128$   &  1.16  & 1.14e+05 & 0.90 & 1.27e+05 & 1.00 & 1.11 & 1.43  \\  
   $256 \times 256$   &  4.51  & 1.14e+05 & 0.90 & 1.27e+05 & 1.00 & 1.11 & 1.42  \\    
  \hline
\end{tabular}
\end{center}
\caption{\textbf{Example 1}.
Minorant $\mathcal{J}^{\ominus}_0$,
majorant $\mathcal{J}^{\oplus}_0$ and
their efficiency indices computed on grids of different sizes.}
\label{tab:Ex1:0}
\end{table}
For $N=3$ 
and $N=8$, the truncation's remainder terms can be
precomputed and are given by 
$\mathcal{E}_3 = 63694.86$
and $\mathcal{E}_8 = 106.06$, respectively.
Since the overall efficiency indices in
Table~\ref{tab:Ex1:global} 
stay all in approximately the same range,
we observe that the method is robust.
However, the efficiency indices for the combined error norm 
$I_{\text{eff}}^{\mathcal{M}_1}$ indicate that the modes $k=1$
and $k=4$ are the most significant to represent the solution by its multiharmonic
approximation.
Comparing the last two 
lines of Table~\ref{tab:Ex1:global} shows  
that the value for representing the cost functional of the exact solution is 
already sufficiently accurate for a truncation index $N=3$.
One of the reasons for this is
that the remainder term $\mathcal{E}_N$ can be precomputed exactly.

\begin{table}[!ht]
\begin{center}
\begin{tabular}{|c|c|cc|cc|c|c|}
  \hline
   grid & $t^{\text{sec}}$ 
   & $\mathcal{J}^{\ominus}_1$ 
   & $I_{\text{eff}}^{\mathcal{J}^\ominus_1}$
   & $\mathcal{J}^{\oplus}_1$ 
   & $I_{\text{eff}}^{\mathcal{J}^\oplus_1}$ 
   & $I_{\text{eff}}^{\mathcal{J},1}$ 
   & $I_{\text{eff}}^{\mathcal{M}_1,1}$ \\
  \hline
   $16   \times   16$   &  0.02  & 4.27e+05 & 0.90 & 4.74e+05 & 1.00 & 1.11 & 6.18  \\ 
   $32   \times   32$   &  0.07  & 4.31e+05 & 0.90 & 4.79e+05 & 1.00 & 1.11 & 6.01  \\ 
   $64   \times   64$   &  0.25  & 4.32e+05 & 0.90 & 4.80e+05 & 1.00 & 1.11 & 5.93  \\   
   $128 \times 128$   &  1.13  & 4.32e+05 & 0.90 & 4.80e+05 & 1.00 & 1.11 & 5.90  \\  
   $256 \times 256$   &  4.66  & 4.32e+05 & 0.90 & 4.80e+05 & 1.00 & 1.11 & 5.89 \\    
  \hline
\end{tabular}
\end{center}
\caption{\textbf{Example 1}.
Minorant $\mathcal{J}^{\ominus}_1$,
majorant $\mathcal{J}^{\oplus}_1$ and
their efficiency indices computed on grids of different sizes.}
\label{tab:Ex1:1}
\end{table}
\begin{table}[!ht]
\begin{center}
\begin{tabular}{|c|c|cc|cc|c|c|}
  \hline
   mode & $t^{\text{sec}}$ 
   & $\mathcal{J}^{\ominus}$ 
   & $I_{\text{eff}}^{\mathcal{J}^\ominus}$
   & $\mathcal{J}^{\oplus}$ 
   & $I_{\text{eff}}^{\mathcal{J}^\oplus}$ 
   & $I_{\text{eff}}^{\mathcal{J}}$ 
   & $I_{\text{eff}}^{\mathcal{M}_1}$ \\
  \hline
   $k=0$   &  4.51  & 1.14e+05 & 0.90 & 1.27e+05 & 1.00 & 1.11 & 1.42  \\   
   $k=1$   &  4.66  & 4.32e+05 & 0.90 & 4.80e+05 & 1.00 & 1.11 & 5.89 \\    
   $k=2$   &  4.75  & 1.79e+05 & 0.90 & 1.99e+05 & 1.00 & 1.11 & 1.64 \\     
   $k=3$   &  4.81  & 6.10e+04 & 0.90 & 6.74e+04 & 1.00 & 1.11 & 1.84 \\    
   $k=4$   &  4.74  & 7.68e+03 & 0.91 & 8.42e+03 & 1.00 & 1.10 & 11.06 \\     
   $k=6$   &  4.72  & 2.05e+02 & 0.90 & 2.29e+02 & 1.00 & 1.11 & 2.08 \\    
   $k=8$   &  4.81  & 1.97e+01 & 0.93 & 2.19e+01 & 1.04 & 1.12 & 1.22 \\       
  \hline
  \hline
   overall ($N=3$) & --  & 2.86e+06 & 0.90 & 3.17e+06 & 1.00 & 1.11 & 2.08 \\   
   overall ($N=8$) & --  & 2.86e+06 & 0.90 & 3.17e+06 & 1.00 & 1.11 & 2.09 \\    
  \hline
\end{tabular}
\end{center}
\caption{\textbf{Example 1}.
Overall minorant $\mathcal{J}^{\ominus}$ and
overall majorant $\mathcal{J}^{\oplus}$, 
their parts, and their efficiency indices computed on a grid of size $256 \times 256$.}
\label{tab:Ex1:global}
\end{table}

\subsubsection{Example 2}
\label{SSSec:Ex2}

We choose the time-analytic, not time-periodic, desired state function
\begin{align}
\label{ydEx2}
 y_d(\boldsymbol{x},t) = e^t \, 0.2  \left((5+ 2 \pi^4) \sin(t) - \cos(t) \right)
			 \sin (x_1 \pi) \sin (x_2 \pi)
\end{align}
having as exact solution the state function
\begin{align}
\label{eqn:exactSol:ex2}
y(\boldsymbol{x},t) = e^t \sin(t) \sin (x_1 \pi) \sin (x_2 \pi).
\end{align}
The approximations by the MhFEM are computed for the truncation index $N = 10$.
For Example 2, it suffices to present here only the overall results as in
Table \ref{tab:Ex1:global}, now presented in 
Table~\ref{tab:Ex2:global}. 
We compare the overall values of majorants and minorants
for different truncation indices $N=6$ 
and $N=10$, for which the corresponding truncation's
remainder terms 
are given by
$\mathcal{E}_6 = 44094.84$ 
and $\mathcal{E}_{10} = 10597.20$.
One can see from the last two 
lines that the truncation index $N = 6$ suffices already to provide an accurate enough approximate solution.
Also efficiency indices being around 1 show that the majorants and minorants perform well for that example.
\begin{table}[!ht]
\begin{center}
\begin{tabular}{|c|c|cc|cc|c|c|}
  \hline
   mode & $t^{\text{sec}}$ 
   & $\mathcal{J}^{\ominus}$ 
   & $I_{\text{eff}}^{\mathcal{J}^\ominus}$
   & $\mathcal{J}^{\oplus}$ 
   & $I_{\text{eff}}^{\mathcal{J}^\oplus}$ 
   & $I_{\text{eff}}^{\mathcal{J}}$ 
   & $I_{\text{eff}}^{\mathcal{M}_1}$ \\
  \hline
   $k=0$   & 4.55 & 3.20e+05     & 0.90 & 3.56e+05     & 1.00 & 1.11 & 1.43 \\    
   $k=1$   & 4.53 & 1.02e+06     & 0.90 & 1.14e+06     & 1.00 & 1.11 & 3.19 \\ 
   $k=2$   & 4.62 & 2.57e+05     & 0.90 & 2.85e+05     & 1.00 & 1.11 & 3.17 \\   
   $k=3$   & 4.80 & 6.07e+04     & 0.91 & 6.69e+04     & 1.00 & 1.10 & 1.82 \\  
   $k=4$   & 4.75 & 1.99e+04     & 0.91 & 2.19e+04     & 1.00 & 1.10 & 1.55 \\   
   $k=6$   & 4.83 & 4.05e+03     & 0.92 & 4.38e+03     & 1.00 & 1.08 & 1.32 \\    
   $k=8$   & 4.90 & 1.30e+03     & 0.93 & 1.40e+03     & 1.00 & 1.07 & 1.20 \\      
   $k=10$ & 4.72 & 5.40e+02     & 0.94 & 5.75e+02     & 1.00 & 1.06 & 1.14 \\  
  \hline
  \hline
   overall ($N=6$)   & --  & 6.34e+06 & 0.90 & 7.06e+06 & 1.00 & 1.11 & 2.09 \\   
   overall ($N=10$) & --  & 6.34e+06 & 0.90 & 7.06e+06 & 1.00 & 1.11 & 2.09 \\  
  \hline
\end{tabular}
\end{center}
\caption{\textbf{Example 2}.
Overall minorant $\mathcal{J}^{\ominus}$ and
overall majorant $\mathcal{J}^{\oplus}$,
their parts, and their efficiency indices computed on a grid of size $256 \times 256$.}
\label{tab:Ex2:global}
\end{table}

\subsubsection{Example 3}
\label{SSSec:Ex3}

The desired state is chosen to be a space-time non-smooth function
\begin{align}
\label{ydEx3}
 y_d(\boldsymbol{x},t) = 
 \chi_{[\frac{1}{2},1]^2}(\boldsymbol{x}) \, \chi_{[\frac{1}{4},\frac{3}{4}]}(t).
\end{align}
Here, we denote by $\chi$ the space-time characteristic function.
The desired state's Fourier cofficients are analytically computed and given by
${y_d^c}_0(\boldsymbol{x}) = \chi_{[\frac{1}{2},1]^2}(\boldsymbol{x})/2$ and
\begin{align}
\label{ydEx3FourierCoeff}
 {y_d^c}_k(\boldsymbol{x})  = \chi_{[\frac{1}{2},1]^2}(\boldsymbol{x}) \frac{
 \left( \sin(\frac{3 k \pi}{2}) -\sin(\frac{k \pi}{2})  \right)}{k \pi} 
 \qquad \text{and} \qquad
 {y_d^s}_k(\boldsymbol{x}) = 0 \qquad
 \forall k \in \mathbb{N}.
\end{align}
The desired state has jumps in space \textbf{and} time.
In this example, the exact solution cannot be precomputed analytically.
Hence, as approximation for it, we use a MhFE representation on the finer grid
with size
$512 \times 512$.
Note also that the Fourier coefficients are zero for the even Fourier modes
besides for $k=0$.
We can observe in Table \ref{tab:Ex3:global} that 
the efficiency indices, especially regarding the
combined norm, are similar for higher modes. We have computed the results
up to $N=11$ 
and also added the results for the modes $k=21$, $k=41$ and $k=81$ to the table
as examples.
The majorants of the combined norm stay approximately in the same range for $k \geq 1$.
The majorants and minorants for the cost functional
are close to 1, 
which demonstrates their efficiency also in this numerical example,
where the given data has jumps in space and time.
\begin{table}[!ht]
\begin{center}
\begin{tabular}{|c|c|cc|cc|c|c|}
  \hline
   mode  & $t^{\text{sec}}$ 
   & $\mathcal{J}^{\ominus}_k$ 
   & $I_{\text{eff}}^{\mathcal{J}^\ominus_k}$
   & $\mathcal{J}^{\oplus}_k$ 
   & $I_{\text{eff}}^{\mathcal{J}^\oplus_k}$ 
   & $I_{\text{eff}}^{\mathcal{J},k}$ 
   & $I_{\text{eff}}^{\mathcal{M}_1,k}$ \\
  \hline
   $k=0$   & 6.44 & 5.71e+04 & 0.92 & 8.20e+04 & 1.32 & 1.44 & 6.48 \\    
   $k=1$   & 6.45 & 9.30e+04 & 0.92 & 1.31e+05 & 1.30 & 1.41 & 3.59 \\  
   $k=3$   & 6.46 & 1.06e+04 & 0.95 & 1.35e+04 & 1.20 & 1.27 & 2.86 \\     
   $k=5$   & 6.39 & 3.90e+03 & 0.97 & 4.63e+03 & 1.15 & 1.19 & 3.16 \\   
   $k=7$   & 6.46 & 2.01e+03 & 0.98 & 2.28e+03 & 1.11 & 1.13 & 3.28 \\     
   $k=9$   & 6.55 & 1.23e+03 & 0.98 & 1.35e+03 & 1.08 & 1.10 & 3.30 \\  
   $k=11$ & 6.48 & 8.23e+02 & 0.99 & 8.89e+02 & 1.07 & 1.08 & 3.23 \\
   \hline      
   $k=21$ & 6.44 & 2.28e+02 & 1.00 & 2.37e+02 & 1.03 & 1.04 & 3.88 \\   
  $k=41$ & 6.46 & 5.99e+01 & 1.00 & 6.19e+01 & 1.03 & 1.03 & 4.81 \\ 
  $k=81$ & 6.42 & 1.53e+01 & 1.00 & 1.63e+01 & 1.06 & 1.06 & 3.19 \\ 
  \hline
\end{tabular}
\end{center}
\caption{\textbf{Example 3}.
Minorants,
majorants,
and their efficiency indices as well as the efficiency indices of the combined norm
computed on a grid of size $256 \times 256$.}
\label{tab:Ex3:global}
\end{table}

\subsection{Numerical results for optimization problem II}
\label{SSec:NumericalResultsII}

We compute the numerical results for the three same cases as for problem I but
now applied on the desired gradient.

\subsubsection{Example 4}
\label{SSSec:Ex4}

We set the desired gradient to be time-periodic and time-analytic
\begin{align*}
 \bsg_d(\boldsymbol{x},t) = \frac{e^t \sin(t) (-3 \cos (t) (\cos (t) +\sin (t))  + (10 \pi^2 + 1 + 2 \pi^4) \sin (t)^2)}{10 \pi}
 \begin{pmatrix} 
 \cos (x_1 \pi) \sin (x_2 \pi)  \\
 \sin (x_1 \pi) \cos (x_2 \pi)
\end{pmatrix}.
\end{align*}
The exact solution for the state function is given by \eqref{eqn:exactSol:ex1}. 
Moreover, we present the efficiency indices for $\tilde{\mathcal{M}}_1^{\oplus}$  
given for the modes by
\begin{align*}
   I_{\text{eff}}^{\tilde{\mathcal{M}}_1,0} =  
  \sqrt{\frac{\tilde{\mathcal{M}}_{1,0}^\oplus(\alpha_0,\beta_0;y_{0 h}^c,p_{0 h}^c,
  \boldsymbol{\tau}_{0 h}^c,\boldsymbol{\rho}_{0 h}^c)}{||||y_{0}^c - y_{0 h}^c||||_{1,0}^2}}
\quad \text{ and } \quad
  I_{\text{eff}}^{\tilde{\mathcal{M}}_1,k} =  
  \sqrt{\frac{\tilde{\mathcal{M}}_{1,k}^\oplus(\alpha_k,\beta_k;\boldsymbol{y}_{k h},\boldsymbol{p}_{k h},
  \boldsymbol{\tau}_{k h},\boldsymbol{\rho}_{k h})}{||||\boldsymbol{y}_{k} - \boldsymbol{y}_{k h}||||_{1,k}^2}}.
\end{align*}
The error norms for the modes are given by
\begin{align*}
||||y_{0}^c - y_{0 h}^c||||_{1,0}^2 =
\left(\frac{1}{2} 
 + \frac{\lambda \underline{\mu_1}^2}{2 C_F^2} \right) \|\nabla y_0^c - \nabla y_{0 h}^c\|_{\Omega}^2
\end{align*}
and
\begin{align*}
||||\boldsymbol{y}_{k} - \boldsymbol{y}_{k h}||||_{1,k}^2 =
\frac{k\omega \lambda \underline{\mu_1}^2}{2 C_F^2} 
  \|\boldsymbol{y}_{k} - \boldsymbol{y}_{k h}\|_{\Omega}^2 +
 \left(\frac{1}{2}  + \frac{\lambda \underline{\mu_1}^2}{2 C_F^2}    \right)
  \|\nabla \boldsymbol{y}_{k} - \nabla \boldsymbol{y}_{k h}\|_{\Omega}^2
\end{align*}
leading the representation
   \begin{align}
  \label{definition:combinedNormH1withktilde}
  \begin{aligned}
 ||||u-v||||_1^2 
  = &\, T \,  ||||y_{0}^c - y_{0 h}^c||||_{1,0}^2 
      + \frac{T}{2} \sum_{k=1}^N  ||||\boldsymbol{y}_{k} - \boldsymbol{y}_{k h}||||_{1,k}^2 
+ \tilde{\mathcal{F}}_N
 \end{aligned}
\end{align}  
with the remainder term
$ \tilde{\mathcal{F}}_N := 
   \frac{T}{2} \sum_{k={N+1}}^\infty  ||||\boldsymbol{y}_{k}||||_{1,k}^2$.
For the numerical experiments, the efficiency index for $\tilde{\mathcal{M}}_1^{\oplus}$ 
from above by estimating \eqref{definition:combinedNormH1withktilde}
from below ignoring the remainder term $ \tilde{\mathcal{F}}_N$ 
leading to the overall efficiency index for $\tilde{\mathcal{M}}_1^{\oplus}$ given by
\begin{align*}
  I_{\text{eff}}^{\tilde{\mathcal{M}}_1} =  
  \sqrt{\frac{
  \tilde{\mathcal{M}}_1^\oplus(\alpha,\beta;\eta,\zeta,\boldsymbol{\tau},\boldsymbol{\rho},v)}{
  T \,  ||||y_{0}^c - y_{0 h}^c||||_{1,0}^2 
      + \frac{T}{2} \sum_{k=1}^N  ||||\boldsymbol{y}_{k} - \boldsymbol{y}_{k h}||||_{1,k}^2}}.
\end{align*}
The corresponding majorants are given by
$\tilde{\mathcal{M}}_{1,0}^\oplus
  = \tilde{\mathcal{M}}_{1,0}^\oplus(\alpha_0,$ 
  $\beta_0;y_{0 h}^c,p_{0 h}^c,\boldsymbol{\tau}_{0 h}^c,\boldsymbol{\rho}_{0 h}^c) 
  = \tilde{\mathcal{J}}_0^\oplus -  \tilde{\mathcal{J}}_0^\ominus 
  + \frac{3 \lambda}{4 C_F^2}(C_F \|{\mathcal{R}_1}^c_0 \|_{\Omega} 
  +  \|{\mathcal{R}_2}^c_0 \|_{\Omega})^2 $
and
$ \tilde{\mathcal{M}}_{1,k}^\oplus 
=  \tilde{\mathcal{M}}_{1,k}^\oplus(\alpha_k,$ $\beta_k;\boldsymbol{y}_{k h},\boldsymbol{p}_{k h},
  \boldsymbol{\tau}_{k h},\boldsymbol{\rho}_{k h})$
$=
   \tilde{\mathcal{J}}_k^\oplus
   -   \tilde{\mathcal{J}}_k^\ominus 
  + \frac{3 \lambda}{4 C_F^2}\left(C_F 
           \|{\mathcal{R}_1}_k\|_{\Omega} + \|{\mathcal{R}_2}_k\|_{\Omega} \right)^2$. 
We present the numerical results for the modes $k=0$ and $k=1$ for different grid sizes
in Tables \ref{tab:Ex4:0} and \ref{tab:Ex4:1}.
The efficiency indices for the majorants
and minorants are very close to 1.00. 
Also the efficiency indices for $\tilde{\mathcal{M}}_{1,0}$ show a good accuracy.
\begin{table}[!ht]
\begin{center}
\begin{tabular}{|c|c|cc|cc|c|c|}
  \hline
   grid & $t^{\text{sec}}$ 
   & $\tilde{\mathcal{J}}^{\ominus}_0$ 
   & $I_{\text{eff}}^{\tilde{\mathcal{J}}^\ominus_0}$
   & $\tilde{\mathcal{J}}^{\oplus}_0$ 
   & $I_{\text{eff}}^{\tilde{\mathcal{J}}^\oplus_0}$ 
   & $I_{\text{eff}}^{\tilde{\mathcal{J}},0}$ 
   & $I_{\text{eff}}^{\tilde{\mathcal{M}}_1,0}$ \\
  \hline
   $16   \times   16$   &  0.02  & 8.92e+03 & 0.99 & 9.85e+03 & 1.09 & 1.10 & 2.04 \\ 
   $32   \times   32$   &  0.06  & 9.24e+03 & 0.99 & 9.95e+03 & 1.07 & 1.08 & 1.97 \\ 
   $64   \times   64$   &  0.24  & 9.32e+03 & 0.99 & 9.97e+03 & 1.05 & 1.07 & 1.94 \\ 
   $128 \times 128$   &  1.04  & 9.34e+03 & 0.98 & 9.98e+03 & 1.05 & 1.07 & 1.93 \\ 
   $256 \times 256$   &  4.39  & 9.35e+03 & 0.98 & 9.98e+03 & 1.05 & 1.07 & 1.92 \\   
  \hline
\end{tabular}
\end{center}
\caption{\textbf{Example 4}.
Minorant $\tilde{\mathcal{J}}^{\ominus}_0$,
majorant $\tilde{\mathcal{J}}^{\oplus}_0$ and
their efficiency indices computed on grids of different sizes.}
\label{tab:Ex4:0}
\end{table}
\begin{table}[!ht]
\begin{center}
\begin{tabular}{|c|c|cc|cc|c|c|}
  \hline
   grid & $t^{\text{sec}}$ 
   & $\tilde{\mathcal{J}}^{\ominus}_1$ 
   & $I_{\text{eff}}^{\tilde{\mathcal{J}}^\ominus_1}$
   & $\tilde{\mathcal{J}}^{\oplus}_1$ 
   & $I_{\text{eff}}^{\tilde{\mathcal{J}}^\oplus_1}$ 
   & $I_{\text{eff}}^{\tilde{\mathcal{J}},1}$ 
   & $I_{\text{eff}}^{\tilde{\mathcal{M}}_1,1}$ \\
  \hline
   $16   \times   16$   & 0.02 & 3.22e+04 & 0.95 & 3.51e+04 & 1.03 & 1.09 & 1.52 \\ 
   $32   \times   32$   & 0.06 & 3.39e+04 & 0.96 & 3.58e+04 & 1.02 & 1.05 & 1.33 \\ 
   $64   \times   64$   & 0.26 & 3.45e+04 & 0.97 & 3.60e+04 & 1.01 & 1.04 & 1.25 \\
   $128 \times 128$   & 1.08 & 3.48e+04 & 0.97 & 3.62e+04 & 1.01 & 1.04 & 1.21 \\
   $256 \times 256$   & 4.31 & 3.51e+04 & 0.98 & 3.65e+04 & 1.01 & 1.04 & 1.19 \\ 
  \hline
\end{tabular}
\end{center}
\caption{\textbf{Example 4}.
Minorant $\tilde{\mathcal{J}}^{\ominus}_1$,
majorant $\tilde{\mathcal{J}}^{\oplus}_1$ and
their efficiency indices computed on grids of different sizes.}
\label{tab:Ex4:1}
\end{table}
Table \ref{tab:Ex4:global} compares
the results for different Fourier modes up
to $N=8$ 
computed on a grid of size $256 \times 256$.
Here, the overall minorants, majorants and efficiency indices are presented,
where
the remainder terms
$\mathcal{E}_N$ for $N=8$ and also $N=6$ have been precomputed exactly.
The values of the efficiency indices vary for different modes $k$.
For example, 
the results for $\tilde{\mathcal{M}}_{1,4}$ indicate that the mode $k=4$ is essential to represent the solution
accurately.
The values for $I_{\text{eff}}^{\mathcal{J}^\ominus_7}$ and $I_{\text{eff}}^{\mathcal{J}^\ominus_8}$
indicate 
that the minorants require a different, higher
refinement for a more accurate representation of the overall solution.
An adaptive scheme is the natural choice. 
On the other hand, 
in these cases the majorants give a good representation for the cost 
functional.
Finally, comparing the last two lines of Table~\ref{tab:Ex4:global} again shows  
that the overall value for representing the cost functional of the exact solution is 
already sufficiently accurate for a truncation index $N=6$.
\begin{table}[!ht]
\begin{center}
\begin{tabular}{|c|c|cc|cc|c|c|}
  \hline
   mode  & $t^{\text{sec}}$ 
   & $\tilde{\mathcal{J}}^{\ominus}$ 
   & $I_{\text{eff}}^{\tilde{\mathcal{J}}^\ominus}$
   & $\tilde{\mathcal{J}}^{\oplus}$ 
   & $I_{\text{eff}}^{\tilde{\mathcal{J}}^\oplus}$ 
   & $I_{\text{eff}}^{\tilde{\mathcal{J}}}$ 
   & $I_{\text{eff}}^{\tilde{\mathcal{M}}_1}$ \\
  \hline
   $k=0$    & 4.39  & 9.35e+03 & 0.98 & 9.98e+03 & 1.05 & 1.07 & 1.92 \\   
   $k=1$    & 4.31  & 3.51e+04 & 0.98 & 3.65e+04 & 1.01 & 1.04 & 1.19 \\  
   $k=2$    & 4.42  & 9.23e+03 & 0.63 & 1.57e+04 & 1.06 & 1.70 & 1.65 \\    
   $k=3$    & 4.43  & 2.85e+03 & 0.58 & 5.06e+03 & 1.03 & 1.78 & 1.08 \\        
   $k=4$    & 4.44  & 5.89e+02 & 0.98 & 6.47e+02 & 1.07 & 1.10 & 5.91 \\      
   $k=6$    & 4.36  & 1.74e+01 & 0.97 & 2.82e+01 & 1.58 & 1.62 & 2.32 \\  
   $k=8$    & 4.33  & 2.99e-01  & 0.23 & 1.37e+00 & 1.05 & 4.59 & 1.97 \\      
  \hline
  \hline
   overall ($N=6$) & --  & 2.09e+05 & 0.89 & 2.45e+05 & 1.04 & 1.17 & 2.46 \\   
   overall ($N=8$) & --  & 2.09e+05 & 0.89 & 2.45e+05 & 1.04 & 1.17 & 2.46 \\    
  \hline
\end{tabular}
\end{center}
\caption{\textbf{Example 4}.
Overall minorant $\tilde{\mathcal{J}}^{\ominus}$ and
overall majorant $\tilde{\mathcal{J}}^{\oplus}$, 
their parts, and their efficiency indices computed on a grid of size $256 \times 256$.}
\label{tab:Ex4:global}
\end{table}

\subsubsection{Example 5}
\label{SSSec:Ex5}

We choose the non time-periodic but time-analytic
desired gradient
\begin{align*}
 \bsg_d(\boldsymbol{x},t) = \frac{-e^t \sin(t) ( 0.1 \cos (t) - \pi^2(1 + 2 \pi^2 0.1))}{\pi}
 \begin{pmatrix} 
 \cos (x_1 \pi) \sin (x_2 \pi)  \\
 \sin (x_1 \pi) \cos (x_2 \pi)
\end{pmatrix}
\end{align*}
leading to
the time-analytic, but not time-periodic 
exact state 
\eqref{eqn:exactSol:ex2}.
We
compute the MhFE approximation
of the desired gradient and solve the systems 
\eqref{equation:MultiFESysBlockk} and \eqref{equation:MultiFESysBlock0}
for 
modes up to $N=10$ 
on a  $256 \times 256$-mesh and present the results 
in Table~\ref{tab:Ex5:global}.
The remainder terms for $N=6$ 
and $N=10$ are 
$\mathcal{E}_6 = 4796.54$ 
and $\mathcal{E}_{10} = 1149.65$, respectively.
The efficiency indices for the overall majorant and minorant show that a truncation index of
$N=6$ already gives a sufficiently accurate approximation for the overall cost functional.
Note that the efficiency index for $\tilde{\mathcal{M}}_{1,2}$ indicates that the mode $k=2$ is essential for
the multiharmonic approximation giving an accurate representation of the solution.
\begin{table}[!ht]
\begin{center}
\begin{tabular}{|c|c|cc|cc|c|c|}
  \hline
   mode  & $t^{\text{sec}}$ 
   & $\tilde{\mathcal{J}}^{\ominus}$ 
   & $I_{\text{eff}}^{\tilde{\mathcal{J}}^\ominus}$
   & $\tilde{\mathcal{J}}^{\oplus}$ 
   & $I_{\text{eff}}^{\tilde{\mathcal{J}}^\oplus}$ 
   & $I_{\text{eff}}^{\tilde{\mathcal{J}}}$ 
   & $I_{\text{eff}}^{\tilde{\mathcal{M}}_1}$ \\
  \hline
   $k=0$   & 4.31 & 2.63e+04 & 1.00 & 2.79e+04 & 1.06 & 1.06 & 1.36 \\   
   $k=1$   & 4.30 & 8.49e+04 & 1.00 & 8.60e+04 & 1.02 & 1.01 & 1.00 \\      
   $k=2$   & 4.39 & 2.08e+04 & 0.98 & 2.21e+04 & 1.04 & 1.06 & 2.83 \\       
   $k=4$   & 4.35 & 1.58e+03 & 0.96 & 1.75e+03 & 1.06 & 1.11 & 1.74 \\      
   $k=6$   & 4.36 & 2.93e+02 & 0.87 & 3.53e+02 & 1.05 & 1.20 & 1.63 \\    
   $k=8$   & 4.37 & 8.95e+01 & 0.82 & 1.20e+02 & 1.10 & 1.34 & 1.08 \\         
   $k=10$ & 4.34 & 3.27e+01 & 0.71 & 5.23e+01 & 1.14 & 1.60 & 1.19 \\   
  \hline
  \hline
   overall ($N=6$)   & --  & 5.23e+05 & 1.00 & 5.43e+05 & 1.04 & 1.04 & 2.00 \\   
   overall ($N=10$) & --  & 5.22e+05 & 1.00 & 5.42e+05 & 1.04 & 1.04 & 2.00 \\  
  \hline
\end{tabular}
\end{center}
\caption{\textbf{Example 5}.
Overall minorant $\tilde{\mathcal{J}}^{\ominus}$ and overall majorant $\tilde{\mathcal{J}}^{\oplus}$,
their parts, and efficiency indices of them and the combined norm computed on a grid of size $256 \times 256$.}
\label{tab:Ex5:global}
\end{table}

\subsubsection{Example 6}
\label{SSSec:Ex6}

We set the space-time non-smooth desired gradient
\begin{align}
\label{ydEx6}
 \bsg_d(\boldsymbol{x},t) = 
 (\chi_{[\frac{1}{2},1]^2}(\boldsymbol{x}) \, \chi_{[\frac{1}{4},\frac{3}{4}]}(t), 
 \chi_{[\frac{1}{2},1]^2}(\boldsymbol{x}) \, \chi_{[\frac{1}{4},\frac{3}{4}]}(t))^T.
\end{align}
Also the coefficients of the Fourier 
expansion associated with $\bsg_d$ 
can be found analytically. They are as in Example 3 given by
\eqref{ydEx3FourierCoeff} for each direction of the gradient \eqref{ydEx6}.
Again the exact solution cannot be 
computed analytically and hence we use its MhFE approximations on the finer mesh of size
$512 \times 512$ 
as a reference.
Table \ref{tab:Ex6:global} presents the results for modes up to truncation index $N=11$ 
as well as for $k=21$, $k=41$ and $k=81$ analogously to Example 3.
The results reflected by the efficiency indices show the good representation by using
the minorants and majorants, especially, considering the efficiency indices in the last two columns
of Table \ref{tab:Ex6:global}.
This again demonstrates the efficiency of the minorants and majorants for data
having jumps in space and time but now for optimization problem II.
\begin{table}[!ht]
\begin{center}
\begin{tabular}{|c|c|cc|cc|c|c|}
  \hline
   mode  & $t^{\text{sec}}$ 
   & $\mathcal{J}^{\ominus}_k$ 
   & $I_{\text{eff}}^{\mathcal{J}^\ominus_k}$
   & $\mathcal{J}^{\oplus}_k$ 
   & $I_{\text{eff}}^{\mathcal{J}^\oplus_k}$ 
   & $I_{\text{eff}}^{\mathcal{J},k}$ 
   & $I_{\text{eff}}^{\mathcal{M}_1,k}$ \\
  \hline
   $k=0$   & 8.12 & 1.55e+01 & 0.79 & 2.05e+01 & 1.04 & 1.32 & 2.23 \\     
   $k=1$   & 8.03 & 7.70e+00 & 0.86 & 1.13e+01 & 1.26 & 1.47 & 2.87 \\  
   $k=3$   & 8.51 & 1.46e+01 & 0.93 & 1.73e+01 & 1.10 & 1.18 & 2.56 \\     
   $k=5$   & 8.25 & 1.44e+01 & 0.99 & 1.65e+01 & 1.14 & 1.15 & 2.59 \\   
   $k=7$   & 8.08 & 8.92e+00 & 0.98 & 9.80e+00 & 1.08 & 1.10 & 1.18 \\     
   $k=9$   & 8.35 & 4.51e+00 & 0.90 & 5.14e+00 & 1.03 & 1.14 & 1.26 \\   
   $k=11$ & 8.36 & 2.57e+00 & 0.96 & 3.13e+00 & 1.17 & 1.22 & 1.51 \\      
    \hline
   $k=21$ & 8.51 & 1.36e+00 & 0.99 & 2.27e+00 & 1.65 & 1.67 & 3.16 \\  
   $k=41$ & 7.78 & 4.54e+00 & 0.86 & 6.26e+00 & 1.18 & 1.38 & 3.09 \\  
   $k=81$ & 7.81 & 3.19e+00 & 0.79 & 7.09e+00 & 1.75 & 2.22 & 3.44 \\  
  \hline
\end{tabular}
\end{center}
\caption{\textbf{Example 6}.
Minorants, majorants, 
and their efficiency indices as well as the efficiency indices of the combined norm
computed on a grid of size $256 \times 256$.}
\label{tab:Ex6:global}
\end{table}

\section{Conclusions and outlook} 
\label{Sec7:Conclusions} 

In this work, the a posteriori error analysis started in
\cite{LangerRepinWolfmayr2016} 
 has been extended now by deriving new  
 lower bounds,
called minorants, for the cost functional leading to an upper estimate for the error
norm of the state and control or equivalently in state and adjoint state.
These lower bounds are guaranteed and computable. Together with using
the results from \cite{Wolfmayr2016} as well as \cite{KollmannKolmbauer:2011}
one can apply the method also to time-periodic optimal control problems, where
box constraints are being imposed on the Fourier coefficients of the control.
The estimates are derived for two different cost functionals,
where the second one is now new in this context.

Since in the linear case the problems are decoupled, the solutions on the Fourier coefficients
could easily be computed on grids of different sizes
depending on the accuracy needed, which could be exactly determined
by using the a posteriori estimates presented in this work
leading to an adaptive method in time. Together with the adaptive finite element method
we then obtain a space-time adaptive method, the adaptive multiharmonic finite element 
method, which we call AMhFEM,
as mentioned for the first time in \cite{LangerRepinWolfmayr2016}.

In this work, a first derivation of preconditioners for the MinRes method for the second optimization
problem has been presented as well as a preconditoner for applying AMLI has been suggested.
Several numerical tests for optimization problem I and II have been presented showing the
efficiency of the upper and especially -- with regard to the article -- lower bounds for the cost functionals in practice.

\section*{Acknowledgments}

The author gratefully acknowledges the financial support by the
Aca\-demy of Finland under the grant 295897.
The author would like to thank the anonymous referees for their
valuable comments 
improving the article.

\bibliographystyle{siam} 
\bibliography{bibliographyarXiv_new3}
\end{document}